\documentclass[twoside,10pt]{article}
\usepackage{graphicx}
\usepackage{subfigure}
\usepackage{url}
\usepackage{color}

\usepackage{amsfonts}
\newcommand\R{\mathbb R}
\newcommand\T{\mathbb T}
\newcommand\Z{\mathbb Z}
\newcommand\ee{\mathrm e}

\newcommand\eps\varepsilon
\newcommand\A{\mathcal A}
\newcommand\F{\mathcal F}
\newcommand\G{\mathcal G}
\newcommand\I{\mathcal I}
\newcommand\J{\mathcal J}
\newcommand\K{\mathcal K}
\newcommand\Lc{\mathcal L}
\newcommand\M{\mathcal M}
\newcommand\Pc{\mathcal P}
\newcommand\Rc{\mathcal R}
\newcommand\Tc{\mathcal T}
\newcommand\W{\mathcal W}
\newcommand\Zc{\mathcal Z}

\newcommand\ds{\displaystyle}
\newcommand\dfrac[2]{\ds\frac{#1}{#2}}
\newcommand\dsqrt[1]{\ds\sqrt{#1}}

\newcommand\p[1]{\left(#1\right)}
\newcommand\pq[1]{\left[#1\right]}
\newcommand\pp[1]{\left\{#1\right\}}
\newcommand\scprod[2]{\left\langle#1,#2\right\rangle}
\newcommand\abs[1]{\left|#1\right|}

\newcommand\tl{\tilde}
\newcommand\wtl{\widetilde}
\newcommand\wh[1]{\widehat{#1}}
\newcommand\ol[1]{\overline{#1}}
\newcommand\df{{\rm d}}
\newcommand\Df{{\rm D}}
\newcommand\deriv[2]{\frac{\df#1}{\df#2}} 
\newcommand\pderiv[2]{\frac{\partial#1}{\partial#2}} 
\newcommand\pdderiv[2]{\frac{\partial^2#1}{\partial#2^{\,2}}}
\newcommand\pdderivm[3]{\frac{\partial^2#1}{\partial#2\,\partial#3}}
\newcommand\Ord{\mathcal O}
\newcommand\rint{\mathop{\rm rint}}
\newcommand\ut{{\rm u}}  
\newcommand\st{{\rm s}}  
\newcommand\loc{{\rm loc}}
\newcommand\vect[2]{\p{\begin{array}{c}#1\\[2pt]#2\end{array}}}
\newcommand\mmatrix[4]{\p{\begin{array}{cc}#1&#2\\[3pt]#3&#4\end{array}}}
\newcommand\symmatrix[3]{\mmatrix{#1}{#2}{#2}{#3}}
\newcommand\tp{^\top}
\newcommand\tr{\mathop{\rm tr}}
\newcommand\ppcf[1]{[\,\ol{#1}\,]}    
\newcommand\epcf[2]{[\,#1,\ol{#2}\,]} 
\newcommand\npcf[1]{[\,#1\,]}         
\newcommand\beps{\bar\eps}    
\newcommand\peps{\bar\eps\,'} 
\newcommand\ceps{\check\eps}  
\newcommand\weps{\wh\eps}     
\newcommand\bg{\ol g}
\newcommand\bS{\ol S}
\newcommand\bh{\ol h}
\newcommand\wtS{\wtl S}
\newcommand\wth{\wtl h}

\newtheorem{theorem}{Theorem}
\newtheorem{numresult}[theorem]{Numerical Result}
\newtheorem{proposition}[theorem]{Proposition}
\newtheorem{lemma}[theorem]{Lemma}
\newcommand\bremark{\noindent{\bf Remark.}\ \ }
\newcommand\eremark{\bigskip}
\newcommand\bremarks{\smallskip\noindent{\bf Remarks.}\ \bnm}
\newcommand\eremarks{\enm\bigskip}
\newcommand\proof{\noindent\textbf{\emph{Proof.}}\quad}
\newcommand\justif{\noindent\textbf{\emph{Justification.}}\quad}
\newcommand\sketchproof{\noindent\textbf{\emph{Sketch of the proof.}}\quad}
\newcommand\proofof[1]{\noindent\textbf{\emph{Proof of #1.}}\quad}
\newcommand\justifof[1]{%
  \noindent\textbf{\emph{Justification of #1.}}\quad}
\newcommand\proofjustifof[2]{%
  \noindent\textbf{\emph{Proof of #1 / Justification of #2.}}\quad}
\newcommand\qed{\ \ \null\nolinebreak\hfill$\frame{\large\phantom a}$}
\newcommand\paragr[1]{\noindent\textbf{\emph{#1.}}\quad}

\newcommand\beq{\begin{equation}}
\newcommand\eeq{\end{equation}}
\newcommand\bea{\begin{eqnarray}}
\newcommand\eea{\end{eqnarray}}
\newcommand\bean{\begin{eqnarray*}}
\newcommand\eean{\end{eqnarray*}}
\newcommand\btm{\vspace{-\baselineskip}\begin{itemize}}
\newcommand\etm{\end{itemize}\vspace{-\baselineskip}}
\newcommand\bnm{\vspace{-\baselineskip}\begin{enumerate}}
\newcommand\enm{\end{enumerate}\vspace{-\baselineskip}}

\begin{document}

\title{Exponentially small splitting of separatrices and transversality
       associated to whiskered tori with quadratic frequency ratio
  \ \footnote{This work has been partially supported by
      the Spanish MINECO-FEDER grant MTM2012-31714,
      the Catalan grant 2014SGR504,
      and the Russian Scientific Foundation grant 14-41-00044
      at the Lobachevsky University of Nizhny Novgorod.
      The author MG has also been supported by the Swedish
      Knut and Alice Wallenberg Foundation grant 2013-0315.
      We also acknowledge the use of EIXAM, the UPC Applied Math cluster system
      for research computing (see \texttt{http://www.ma1.upc.edu/eixam/}),
      and in particular Pau Rold\'an and Albert Granados for their support
      in the use of this cluster.}}
\author{\sc
    Amadeu Delshams$\,^1$,
  \ Marina Gonchenko$\,^2$,
\\[4pt]\sc
  \ Pere Guti\'errez$\,^1$
\\[12pt]
  {\small
  $^1\;$\parbox[t]{5.2cm}{
    Dep. de Matem\`atiques\\
    Universitat Polit\`ecnica de Catalunya\\
    Av. Diagonal 647, 08028 Barcelona\\
    {\footnotesize
      \texttt{amadeu.delshams@upc.edu}\\
      \texttt{pere.gutierrez@upc.edu}}}
  \qquad
  $^2\;$\parbox[t]{5.2cm}{
    Matematiska institutionen\\
    Uppsala universitet\\
    box 480, 751 06 Uppsala\\
    {\footnotesize
      \texttt{marina.gonchenko@math.uu.se}}}
  }}
\maketitle
\begin{abstract}
The splitting of invariant manifolds of whiskered (hyperbolic) tori with two
frequencies in a nearly-integrable Hamiltonian system, whose hyperbolic part is
given by a pendulum, is studied. We consider a torus with a fast frequency
vector $\omega/\sqrt\varepsilon$, with $\omega=(1,\Omega)$ where the frequency
ratio $\Omega$ is a quadratic irrational number. Applying the
Poincar\'e-Melnikov method, we carry out a careful study of the dominant
harmonics of the Melnikov potential. This allows us to provide an asymptotic
estimate for the maximal splitting distance, and show the existence of
transverse homoclinic orbits to the whiskered tori with an asymptotic estimate
for the transversality of the splitting. Both estimates are exponentially small
in $\varepsilon$, with the functions in the exponents being periodic with
respect to $\ln\varepsilon$, and can be explicitly constructed from the
continued fraction of $\Omega$. In this way, we emphasize the strong dependence
of our results on the arithmetic properties of $\Omega$. In particular, for
quadratic ratios $\Omega$ with a 1-periodic or 2-periodic continued fraction
(called metallic and metallic-colored ratios respectively), we provide accurate
upper and lower bounds for the splitting. The estimate for the maximal
splitting distance is valid for all sufficiently small values of $\varepsilon$,
and the transversality can be established for a majority of values of
$\varepsilon$, excluding small intervals around some transition values where
changes in the dominance of the harmonics take place, and bifurcations could
occur.
\par\vspace{12pt}
\noindent\emph{Keywords}:
  splitting of separatrices,
  transverse homoclinic orbits,
  Melnikov integrals,
  quadratic frequency ratio.
\end{abstract}

\section{Introduction and setup}

\subsection{Background and state of the art}

This paper is dedicated to the study of the exponentially small
splitting of separatrices in a perturbed
3-degree-of-freedom Hamiltonian system,
associated to a 2-dimensional whiskered torus
(invariant hyperbolic torus) whose frequency ratio
is an arbitrary \emph{quadratic irrational number}
(i.e.~a real root of a quadratic polynomial with integer coefficients).

We start with an integrable Hamiltonian $H_0$
having whiskered (hyperbolic) tori with coincident
stable and unstable whiskers (invariant manifolds).
We focus our attention on a torus,
with a frequency vector of \emph{fast frequencies}:
\beq
  \omega_\eps = \frac\omega{\sqrt{\eps}}\;,
  \qquad
  \omega=(1,\Omega),
\label{eq:omega_eps}
\eeq
whose \emph{frequency ratio} $\Omega$ is a quadratic irrational number.
If we consider a perturbed Hamiltonian $H=H_0+\mu H_1$, where $\mu$ is small,
in general the whiskers do not coincide anymore.
This phenomenon has got the name of \emph{splitting of separatrices},
and is related to the non-integrability of the system
and the existence of chaotic dynamics.
If we assume, for the two involved parameters,
a relation of the form $\mu=\eps^r$
for some $r>0$, we have a problem of singular perturbation
and in this case the splitting is \emph{exponentially small}
with respect to $\eps$.
Our aim is to detect homoclinic orbits
(i.e.~intersections between the stable and unstable manifolds)
associated to persistent whiskered tori,
provide an \emph{asymptotic estimate} for both the
\emph{maximal splitting distance}
and its \emph{transversality}, and show the dependence of such estimates on the
\emph{arithmetic properties} of the frequency ratio~$\Omega$.

To measure the splitting, it is very usual to apply
the \emph{Poincar\'e--Melnikov method}, introduced by
Poincar\'e in \cite{Poincare90} and rediscovered much later by
Melnikov and Arnold \cite{Melnikov63,Arnold64}.
By considering a transverse section to the
stable and unstable perturbed whiskers,
one can consider a function $\M(\theta)$, $\theta\in\T^2$,
usually called \emph{splitting function},
giving the vector distance, with values in $\R^2$,
between the whiskers on this section,
along the complementary directions.
The method provides a first order approximation
to this function, with respect to the parameter $\mu$,
given by the (vector) \emph{Melnikov function} $M(\theta)$,
defined by an integral (see for instance \cite{Treschev94}).
We have
\beq\label{eq:melniapproxM}
  \M(\theta)=\mu M(\theta)+\Ord(\mu^2),
\eeq
and hence for $\mu$ small enough the simple zeros of $M(\theta)$
give rise to simple zeros of $\M(\theta)$,
i.e.~to transverse intersections between the perturbed whiskers.
In this way, we can obtain asymptotic estimates for both the
maximal splitting distance as the maximum
of the function $\abs{\M(\theta)}$,
and for the transversality of the splitting,
which can be measured by the minimal eigenvalue (in modulus)
of the $(2\times2)$-matrix $\Df\M(\theta^*)$, for any given zero $\theta^*$.

An important fact, related to the Hamiltonian character of the system,
is that both functions $\M(\theta)$ and $M(\theta)$
are gradients of scalar functions \cite{Eliasson94,DelshamsG00}:
\beq\label{eq:gradients}
  \M(\theta)=\nabla\Lc(\theta),
  \qquad
  M(\theta)=\nabla L(\theta).
\eeq
Such scalar funtions are called \emph{splitting potential} and
\emph{Melnikov potential} respectively.
This means that there always exist homoclinic orbits, which correspond to
critical points of $\Lc$.

As said before, the case of fast frequencies $\omega_\eps$
as in~(\ref{eq:omega_eps}), with a perturbation of order $\mu=\eps^r$,
for a given $r$ as small as possible,
turns out to be a \emph{singular problem}.
The difficulty comes from the fact that the Melnikov function $M(\theta)$
is exponentially small in $\eps$,
and the Poincar\'e--Melnikov method can be directly applied
only if one assumes that $\mu$ is exponentially small with respect to $\eps$.
In order to validate the method in the case $\mu=\eps^r$,
one has to ensure that the error term is also exponentially small,
and the Poincar\'e--Melnikov approximation dominates it.
To overcome such a~difficulty in the study of the exponentially small splitting,
it was introduced in \cite{Lazutkin03}
the use of parameterizations of a complex strip of the whiskers
(whose width is defined by the singularities of the unperturbed ones)
by periodic analytic functions, together with flow-box coordinates.
This tool was initially developed for the Chirikov standard map
\cite{Lazutkin03}, and allowed several authors
to validate the Poincar\'e--Melnikov method
for Hamiltonians with one and a~half degrees of freedom (with 1~frequency)
\cite{HolmesMS88,Scheurle89,DelshamsS92,DelshamsS97,Gelfreich97}
and for area-preserving maps \cite{DelshamsR98}.

Later, those methods were extended to the case of whiskered tori
with 2~frequencies. In this case, the arithmetic properties of the frequencies
play an important role in the exponentially small asymptotic estimates of the
splitting function, due to the presence of \emph{small divisors}.
This was first mentioned in \cite{Lochak90}, later
detected in \cite{Simo94}, and rigorously proved in
\cite{DelshamsGJS97} for the quasi-periodically forced pendulum,
assuming a polynomial perturbation in the coordinates associated
to the pendulum. Recently, a more general (meromorphic) perturbation
has been considered in \cite{GuardiaS12}.
It is worth mentioning that,
in some cases, the Poincar\'e--Melnikov method does not predict correctly
the size of the splitting, as shown in \cite{BaldomaFGS12},
where a Hamilton--Jacobi method is instead used.
This method was previously used in \cite{Sauzin01,LochakMS03,RudnevW00},
where exponentially small estimates for the transversality of the splitting
were obtained, excluding
some intervals of the perturbation parameter~$\eps$.
Similar results were obtained in \cite{DelshamsG04,DelshamsG03}.
Moreover, the \emph{continuation} of the transversality for all
sufficiently small values of $\eps$ was shown in \cite{DelshamsG04}
for the concrete case of
the famous \emph{golden ratio} $\Omega=(\sqrt5-1)/2$,
and in \cite{DelshamsGG14c} for the case of
the \emph{silver ratio} $\Omega=\sqrt2-1$,
provided certain conditions on the phases of the perturbation are fulfilled.
Otherwise, homoclinic bifurcations can occur, as studied,
for instance, in \cite{SimoV01} for the Arnold's example.
Let us also mention that analogous estimates could also be obtained
from a careful averaging out of the fast angular variables \cite{ProninT00},
at least concerning sharp upper bounds of the splitting.

In general, in the quoted papers the frequency ratio is assumed to be
a given concrete quadratic number (golden, silver).
A~generalization to some other concrete quadratic frequency ratios
was considered in \cite{DelshamsG03},
extending the asymptotic estimates for the maximal
splitting distance, but without a satisfatory result concerning transversality.
Recently, a parallel study of the cases of~2 and 3~frequencies has been
considered in \cite{DelshamsGG14a}
(in the case of 3~frequencies, with a frequency vector
$\omega=(1,\Omega,\Omega^2)$, where $\Omega$ is a concrete
cubic irrational number),
obtaining also exponentially small estimates
for the maximal splitting distance.

In this paper, we consider a 2-dimensional torus whose frequency
ratio $\Omega$ in~(\ref{eq:omega_eps}) is a given quadratic irrational number.
Our main objective is to develop a methodology,
allowing us to obtain asymptotic estimates for both
the maximal splitting distance and the transversality of the splitting,
whose dependence on $\eps$ is described
by two \emph{piecewise-smooth functions}
denoted $h_1(\eps)$ and $h_2(\eps)$, respectively (see Theorem~\ref{thm:main}).
Such functions are both periodic with respect to $\ln\eps$,
and their behavior depends strongly on the arithmetic properties
of the frequency ratio $\Omega$.
In particular, we show that the functions $h_1(\eps)$ and $h_2(\eps)$
can be constructed explicitly from the \emph{continued fraction} of~$\Omega$,
and we can deduce some of their
properties like the number of corners in each period
(this can be seen as an indication of the complexity of the dependence
on $\eps$ of the splitting).
Our goal is to show that our methods can be applied to an arbitrary
quadratic ratio, and hence the results
on the splitting distance and transversality generalize the ones of
\cite{DelshamsG03,DelshamsG04,DelshamsGG14a,DelshamsGG14c}.
Although we do not study here the continuation
of the transversality for all values of $\eps\to0$,
we stress that this could be carried out
by means of a specific study in each case,
as done in \cite{SimoV01,DelshamsG04,DelshamsGG14c} for some concrete
(golden, silver) ratios.

We point out that the periodicity in $\ln\eps$ of the functions
$h_1(\eps)$ and $h_2(\eps)$ comes directly from the special properties
of the continued fractions of quadratic numbers
and cannot be satisfied in other cases
(see \cite{DelshamsGG14b}, where the case of numbers of constant type
is considered).

As we describe in Section~\ref{sect:main},
the dependence on $\eps$ of the function $h_2(\eps)$, associated to
the asymptotic estimates for the transversality of the splitting,
is more cumbersome than for the function $h_1(\eps)$, associated to
the maximal splitting distance.
We stress here that, for some purposes, it is not necessary
to establish the transversality of the splitting,
and it can be enough to provide estimates of the maximal splitting distance.
Indeed, such estimates imply the existence of splitting between the invariant
manifolds, which provides a strong indication of the non-integrability
of the system near the given torus, and opens the door to the application
of topological methods \cite{GideaR03,GideaL06} for the study of
Arnold diffusion in such systems.

\subsection{Setup}

Here we describe the nearly-integrable
Hamiltonian system under consideration.
In particular, we study a \emph{singular} or
\emph{weakly hyperbolic} (\emph{a priori stable})
Hamiltonian with 3 degrees of freedom possessing a
2-dimensional whiskered torus with fast frequencies.
In canonical coordinates
$(x,y, \varphi, I)\in\T\times \R \times \T^2
\times \R^2$, with the symplectic form
$\df x \wedge\df y +\df\varphi\wedge\df I$, the Hamiltonian is defined by
\bea
&&H(x,y, \varphi, I) = H_0 (x,y, I) + \mu H_1(x, \varphi),
\label{eq:HamiltH}
\\
&&H_0 (x, y, I) =
\langle \omega_\eps, I\rangle + \frac{1}{2} \langle\Lambda I, I\rangle
+\frac{y^2}{2} + \cos x -1,
\label{eq:HamiltH0}
\\
&&H_1 (x, \varphi)= h(x) f(\varphi).
\label{eq:HamiltH1}
\eea
Our system has two parameters $\eps>0$ and $\mu$, linked
by a relation $\mu=\eps^r$, $r>0$ (the smaller $r$ the better).
Thus, if we consider $\eps$ as the unique parameter, we have
a singular problem for $\eps\to 0$.
See \cite{DelshamsG01} for a discussion about singular and regular problems.

Recall that we are assuming a vector of fast frequencies
$\omega_\eps = \omega/\sqrt{\eps}$ as given
in~(\ref{eq:omega_eps}), with a frequency vector $\omega=(1,\Omega)$
whose frequency ratio $\Omega$ is a \emph{quadratic irrational number};
we assume without loss of generality that $0<\Omega<1$.
It is a well-known property (and we prove it in Section~\ref{sect:primary};
see also \cite[\S II.2]{Lang95})
that any vector with quadratic ratio satisfies a \emph{Diophantine condition}
\beq\label{eq:DiophCond}
|\langle k, \omega\rangle| \geq \frac{\gamma}{\abs k}\,,
\qquad \forall k\in \Z^2\setminus\pp0,
\eeq
with some $\gamma>0$.
We also assume in~(\ref{eq:HamiltH0}) that $\Lambda$ is a symmetric
($2\times2$)-matrix, such that $H_0$ satisfies the condition of
\emph{isoenergetic nondegeneracy}
\beq
    \det \left(
    \begin{array}{cc}
    \Lambda & \omega\\
    \omega\tp & 0
    \end{array}
    \right) \neq 0.
\label{eq:isoenerg}
\eeq

For the perturbation $H_1$ in~(\ref{eq:HamiltH1}),
we deal with the following analytic periodic functions,
\beq\label{eq:hf}
h(x) = \cos x,
\qquad
f(\varphi)= \sum_{k\in\Zc}
e^{-\rho\abs k} \cos(\langle k, \varphi \rangle - \sigma_k),
\quad\mbox{with }\ \sigma_k\in\T,
\eeq
where we introduce, in order to avoid repetitions in the Fourier series,
the set
\beq\label{eq:calZ}
  \Zc = \{k=(k_1, k_2)\in \Z^2: k_2>0 \mbox{ or }(k_2=0, k_1\ge 0)\}.
\eeq
The constant $\rho>0$ gives the complex width of analyticity
of the function $f(\varphi)$.
Concerning the phases $\sigma_k$, they can be chosen arbitrarily
for the purpose of this paper.

To justify the form of the perturbation $H_1$
chosen in~(\ref{eq:HamiltH1}) and~(\ref{eq:hf}),
we stress that it makes easier the explicit
computation of the Melnikov potential, which
is necessary in order to show that it dominates
the error term in~(\ref{eq:melniapproxM}),
and therefore to establish the existence of splitting.
Moreover, the fact that all coefficients \ $f_k=e^{-\rho\abs k}$,
\ in the Fourier expansion with respect to $\varphi$,
are nonzero and have an exponential decay, ensures
that the study of the dominant harmonics of the Melnikov potential
can be carried out directly from the arithmetic properties
of the frequency vector~$\omega$ (see Section~\ref{sect:asympt_est}).
Since our method is completely constructive, a perturbation with
another kind of concrete harmonics $f_k$ could also be considered
(like $f_k=\abs k^me^{-\rho\abs k}$), simply at the cost of more cumbersome
computations in order to determine the dominant harmonics
of the Melnikov potential.

It is worth reminding that the Hamiltonian
defined in~(\ref{eq:HamiltH}--\ref{eq:hf}) is paradigmatic,
since it is a generalization of the famous Arnold's example
(introduced in \cite{Arnold64} to illustrate the transition
chain mechanism in Arnold diffusion).
It provides a model for the behavior of
a nearly-integrable Hamiltonian system near a single resonance
(see \cite{DelshamsG01} for a motivation),
and has often been considered in the literature
(see for instance \cite{GallavottiGM99b,LochakMS03,DelshamsGS04}).
We also mention that a perturbation with an exponential decay as
the function $f(\varphi)$ in~(\ref{eq:hf}) has also been considered
(see for instance \cite{ProninT00}).
In the present paper, our aim is to emphasize
the dependence of the splitting, and
its transversality, on the arithmetic properties
of the frequency vector~$\omega$.

Let us describe the invariant tori and whiskers,
as well as the splitting and Melnikov functions.
First, it is clear that the unperturbed system $H_0$
(that corresponds to $\mu=0$)
consists of the pendulum given by $P(x,y)=y^2/2+\cos x-1$, and
2~rotors with fast frequencies:
\ $\dot{\varphi}= \omega_\eps + \Lambda I$, \ $\dot{I}=0$.
\ The pendulum has a hyperbolic equilibrium at the origin, with
separatrices that correspond to the curves given by $P(x,y)= 0$.
We parameterize the upper separatrix of the pendulum as
$(x_0(s), y_0(s))$, $s\in \R$, where
\[
  x_0(s)=4 \arctan e^s, \qquad y_0(s) = \frac{2}{\cosh s}\,.
\]
Then, the lower separatrix has the parametrization $(x_0(-s), -y_0(-s))$.
For the rotors system $(\varphi, I)$, the solutions are
\ $I= I_0$, \ $\varphi = \varphi_0+t(\omega_\eps + \Lambda I_0)$.
\ Consequently, the Hamiltonian $H_0$ has a 2-parameter family of 2-dimensional
whiskered tori: in coordinates $(x,y,\varphi,I)$,
each torus can be parameterized as
\[
  \Tc_{I_0}:\qquad(0,0,\theta,I_0),\quad\theta\in\T^2,
\]
and the inner dynamics on each torus is
\ $\dot\theta = \omega_\eps+\Lambda I_0$.
\ Each invariant torus has a \emph{homoclinic whisker},
i.e.~coincident 3-dimensional stable and unstable invariant manifolds,
which can be parameterized as
\beq\label{eq:W0}
\W_{I_0}:\qquad(x_0(s), y_0(s),\theta,I_0),
\quad s\in \R, \ \theta\in\T^2,
\eeq
with the inner dynamics given by
\ $\dot{s}=1$, \ $\dot{\theta}=\omega_\eps+\Lambda I_0$.

In fact, the collection of the whiskered tori
for all values of $I_0$ is
a 4-dimensional \emph{normally hyperbolic invariant manifold},
parameterized by $(\theta,I)\in\T^2\times\R^2$.
This manifold has a 5-dimensional homoclinic manifold,
which can be parameterized by $(s,\theta,I)$,
with inner dynamics 
$\dot s=1$, \ $\dot\theta=\omega_\eps+\Lambda I$, \ $\dot I=0$.
\ We stress that this approach is usually considered
in the study of Arnold diffusion
(see for instance \cite{DelshamsLS06}).

Among the family of whiskered tori and homoclinic whiskers,
we will focus our attention on the torus $\Tc_0$,
whose frequency vector is $\omega_\eps$ as in~(\ref{eq:omega_eps}),
and its associated homoclinic whisker $\W_0$.

When adding the perturbation $\mu H_1$, for $\mu\ne0$ small enough
the \emph{hyperbolic KAM theorem} can be applied
(see for instance \cite{Niederman00})
thanks to the Diophantine condition~(\ref{eq:DiophCond})
and to the isoenergetic nondegeneracy~(\ref{eq:isoenerg}).
For $\mu$ small enough, the whiskered torus persists
with some shift and deformation, as a perturbed torus $\Tc=\Tc^{(\mu)}$,
as well as its local whiskers $\W_\loc=\W^{(\mu)}_\loc$
(a precise statement can be found, for instance,
in \cite[Th.~1]{DelshamsGS04}).

The local whiskers can be extended along the flow,
but in general for $\mu\ne0$ the \emph{(global) whiskers}
do not coincide anymore, and one expects the existence of splitting between
the (3-dimensional) stable and unstable whiskers,
denoted $\W^\st=\W^{\st,(\mu)}$ and $\W^\ut=\W^{\ut,(\mu)}$ respectively.
Using \emph{flow-box coordinates}
(see \cite{DelshamsGS04}, where the $n$-dimensional case is considered)
in a neighbourhood containing a piece of
both whiskers (away from the invariant torus),
one can introduce parameterizations of the perturbed whiskers,
with parameters $(s,\theta)$ inherited
from the unperturbed whisker~(\ref{eq:W0}),
and the inner dynamics
\[
  \dot s=1, \qquad \dot{\theta}=\omega_\eps.
\]
Then, the distance between the stable whisker $\W^\st$ and
the unstable whisker $\W^\ut$ can be measured by comparing
such parameterizations along the complementary directions.
The number of such directions is~3 but,
due to the energy conservation, it is enough to consider 2~directions,
say the ones related to the action coordinates.
In this way, one can introduce a (vector) \emph{splitting function},
with values in $\R^2$, as the difference of the parameterizations
$\J^{\st,\ut}(s,\theta)$ of (the action components of)
the perturbed whiskers $\W^\st$ and $\W^\ut$.
Initially this splitting function depends on $(s,\theta)$,
but it can be restricted to a transverse section
by considering a fixed $s$, say $s=0$,
and we can define as in \cite[\S5.2]{DelshamsG00} the splitting function
\beq\label{eq:defM}
  \M(\theta):=\J^\ut(0,\theta)-\J^\st(0,\theta),
  \quad\theta\in\T^2.
\eeq
As said in~(\ref{eq:gradients}),
this function turns out to be the gradient
of the (scalar) \emph{splitting potential} $\Lc(\theta)$
(see \cite{DelshamsG00,Eliasson94}).
Notice that the \emph{nondegenerate critical points} of~$\Lc$ correspond
to simple zeros of~$\M$ and give rise
to \emph{transverse homoclinic orbits} to the whiskered torus.

Applying the Poincar\'e--Melnikov method,
the first order approximation~(\ref{eq:melniapproxM}) of the
splitting function is given by
the (vector) \emph{Melnikov function} $M(\theta)$,
which is the gradient of the \emph{Melnikov potential}:
\ $M(\theta)=\nabla L(\theta)$.
\ The latter one can be defined as an integral:
we consider any homoclinic trajectory of
the unperturbed homoclinic whisker $\W_0$ in~(\ref{eq:W0}),
starting on the section $s=0$, and the trajectory on the torus $\Tc_0$
to which it is asymptotic as $t\to\pm\infty$,
and we substract the values of the perturbation~$H_1$ on the two trajectories.
This gives an absolutely convergent integral,
depending on the initial phase~$\theta\in\T^2$ of the considered trajectories:
\bea
  \nonumber
  L(\theta)
  &:=
  &-\int_{-\infty}^{\infty}
    [H_1(x_0(t),\theta+t\omega_\eps)-H_1(0,\theta+t\omega_\eps)]\,\df t
\\
  \label{eq:L}
  &=
  &-\int_{-\infty}^{\infty}
    [h(x_0(t))-h(0)] f(\theta+t\omega_\eps)\,\df t,
\eea
where we have taken into account the specific form~(\ref{eq:HamiltH1})
of the perturbation.

Our choice of the pendulum $P(x,y)=y^2/2+\cos x-1$ in~(\ref{eq:HamiltH0}),
whose separatrix has simple poles,
makes it possible to use the method of residues in order to compute
the coefficients $L_k$ of the Fourier expansion of
the Melnikov potential~$L(\theta)$.
Such coefficients turn out to be exponentially small in $\eps$
(see their expression in Section~\ref{sect:gn}).
For each value of~$\eps$ only a finite
number of \emph{dominant harmonics} are relevant
to find the nondegenerate critical points of $L(\theta)$,
i.e.~the simple zeros of the Melnikov function $M(\theta)$.
Due to the exponential decay of the Fourier coefficients of $f(\varphi)$
in~(\ref{eq:hf}),
it is not hard to study such dominance and its dependence on $\eps$.

In order to give asymptotic estimates for both the maximal splitting distance
and the transversality of the homoclinic orbits,
the estimates obtained for the Melnikov function $M(\theta)$
have to be validated also for the splitting function~$\M(\theta)$.
The difficulty in the application of 
the Poincar\'e--Melnikov approximation~(\ref{eq:melniapproxM}),
due to the exponential smallness in $\eps$
of the function $M(\theta)$ in our singular case $\mu=\eps^r$,
can be solved by obtaining upper bounds (on a complex domain)
for the \emph{error term} in~(\ref{eq:melniapproxM}),
showing that, if $r>r^*$ with a suitable $r^*$, its Fourier coefficients
are dominated by the coefficients of $M(\theta)$
(see also~\cite{DelshamsGS04}).

We stress that our approach can also be directly applied to other
classical 1-degree-of-freedom Hamiltonians $P(x,y)=y^2/2+V(x)$,
with a potential $V(x)$ having a unique nondegenerate maximum,
although the use of residues becomes more cumbersome when the
complex parameterization of the separatrix has poles of higher orders
(see some examples in \cite{DelshamsS97}).

\subsection{Main result}\label{sect:main}

For the Hamiltonian system~\mbox{(\ref{eq:HamiltH}--\ref{eq:hf})}
with the 2 parameters linked by
$\mu=\eps^r$, $r>r^*$ (with some suitable $r^*$),
and a frequency vector $\omega=(1,\Omega)$ with a quadratic ratio $\Omega$,
our main result provides exponentially small \emph{asymptotic estimates}
for some measures of the splitting.
On one hand, we obtain an asymptotic estimate
for the \emph{maximal distance} of splitting,
given in terms of the maximum size in modulus of the
splitting function $\M(\theta)=\nabla\Lc(\theta)$,
and this estimate is valid for all $\eps$ sufficiently small.
On the other hand, we show the existence
of \emph{transverse homoclinic orbits},
given as simple zeros $\theta^*$ of $\M(\theta)$
(or, equivalently, as nondegenerate critical points of $\Lc(\theta)$),
and we obtain an asymptotic estimate for the \emph{transversality} of the
homoclinic orbits, measured by the minimal eigenvalue (in modulus)
of the matrix $\Df\M(\theta^*)=\Df^2\Lc(\theta^*)$,
at each zero of $\M(\theta)$.
This result on transversality is valid for \emph{``almost all''} $\eps$
sufficiently small,
since we have to exclude a small neighborhood of
a finite number of geometric sequences
where homoclinic bifurcations could take place.

With our approach, the Poincar\'e--Melnikov method can be validated
for an exponent $r>r^*$ with $r^*=3$,
although a lower value of $r^*$ can be given in some particular cases
(see remark~\ref{rk:hf2} after Theorem~\ref{thm:main}).
However, such values of $r^*$ are not optimal and could be improved
using other methods, like the parametrization of the
whiskers as solutions of Hamilton--Jacobi equation
(see for instance \cite{LochakMS03,BaldomaFGS12}).
In this paper, the emphasis is put in the generalization of the estimates to
the case of an arbitrary quadratic frequency ratio $\Omega$,
rather than in the improvement of the value of $r^*$.

Due to the form of $f(\varphi)$ in~(\ref{eq:hf}),
the Melnikov potential $L(\theta)$
is readily presented in its Fourier series (see Section~\ref{sect:gn}),
with coefficients $L_k=L_k(\eps)$
which are exponentially small in $\eps$.
We use this expansion of $L(\theta)$ in order to detect its
\emph{dominant harmonics} for every given $\eps$.
A careful control of the error term in~(\ref{eq:melniapproxM})
ensures that the dominant harmonics of $L(\theta)$ correspond
to the dominant harmonics of the splitting potential $\Lc(\theta)$.
Such a dominance is also valid for the splitting function
$\M(\theta)$, since the size of their Fourier coefficients
$\M_k$ (vector) and $\Lc_k$ (scalar) is directly related:
\ $\abs{\M_k}=\abs k\,\Lc_k$, \ $k\in\Zc$
(recall the definition~(\ref{eq:calZ})).

As shown in Section~\ref{sect:technical}, in order to obtain an
asymptotic estimate for the maximal distance of splitting,
it is enough to consider the first dominant harmonic,
given by some vector in $\Zc$ which depends on the perturbation parameter:
$k=S_1(\eps)$.
Using estimates for this dominant harmonic $\Lc_{S_1}$
as well as for all the remaining harmonics,
we show that the dominant harmonic
is large enough to ensure that it provides an approximation
to the maximum size of the whole splitting function
(see also \cite{DelshamsGG14a,DelshamsGG14b}).
On the other hand, to show the transversality of the splitting,
it is necessary to consider at least two dominant harmonics in order to
prove the nondegeneracy of the critical points of the splitting potential
(see also \cite{DelshamsG03,DelshamsG04}).
For most values of the parameter $\eps$, it is enough to consider
the two \emph{``essential''} dominant harmonics $\Lc_{S_1}$ and $\Lc_{S_2}$
(i.e.~the two most dominant harmonics whose
associated vectors $S_1(\eps),S_2(\eps)\in\Zc$
are linearly independent, see Section~\ref{sect:resonantseq}),
and the transversality is then directly established.

However, one has to consider at least three harmonics
for $\eps$ near to some \emph{``transition values''} $\weps$,
introduced below as the values
at which a change in the \emph{second} essential dominant harmonic occurs and,
consequently, the second and some subsequent harmonics have similar sizes.
Such transition values turn out to be \emph{corners}
of the function $h_2(\eps)$,
related to the size of the second dominant harmonic (see the theorem below),
and are given by a finite number of geometric sequences.
The study of the transversality for $\eps$ close to a transition value,
which is not considered in this paper, requires to carry out a specific study
that depends strongly on the quadratic frequency ratio~$\Omega$,
and on the concrete perturbation considered in~(\ref{eq:hf}).
In some cases, the \emph{continuation} of the transversality
for all sufficiently small values $\eps\to0$ can be established
under a certain condition on the phases $\sigma_k$ in~(\ref{eq:hf}),
as done in \cite{DelshamsG04} and \cite{DelshamsGG14c}
for the golden and silver ratios, respectively.
Otherwise, \emph{homoclinic bifurcations} can occur when $\eps$ goes accross
a~transition value (see for instance \cite{SimoV01}).

The dependence on $\eps$ of the size of the splitting and its transversality
is closely related to the arithmetic properties of the
frequency vector $\omega=(1,\Omega)$, since
the integer vectors $k\in\Zc$ associated to the dominant harmonics
can be found, for any $\eps$, among
the main quasi-resonances of the vector $\omega$,
i.e.~the vectors $k$ giving the
``least'' small divisors $\abs{\scprod k\omega}$
(relatively to the size of $\abs k$).
In Section~\ref{sect:quadrfreq}, we develop
a methodology for a complete study of the \emph{resonant properties}
of vectors with a quadratic ratio,
which is one of the main goals of this paper.
This methodology relies in the classification,
established in \cite{DelshamsG03} for any vector
with a quadratic ratio $\Omega$,
of the integer vectors~$k$ into \emph{``resonant sequences''}
(see also Sections~\ref{sect:resonantseq} and~\ref{sect:primary}
for definitions).
Among them, the sequence of \emph{primary resonances}
corresponds to the vectors $k$ which fit best
the Diophantine condition~(\ref{eq:DiophCond}),
and the vectors $k$ belonging to the remaining sequences
are called \emph{secondary resonances}.

As particular cases, for the golden ratio $\Omega=(\sqrt5-1)/2$
the primary resonances can be described in terms of the Fibonacci numbers:
$k=(-F_{n-1},F_n)$ (see for instance \cite{DelshamsG04}),
and in the case of the silver ratio $\Omega=\sqrt2-1$
the primary resonances are given in terms of the Pell numbers
(see \cite{DelshamsGG14c}), which play an analogous role.
In general, for a given quadratic ratio $\Omega$
the sequence of primary resonances,
as well as the remaining resonant sequences,
can be determined from the continued fraction of $\Omega$,
which is eventually periodic, i.e.~periodic starting at some element
(see Section~\ref{sect:contfract}).
We can construct, from the continued fraction, a \emph{unimodular} matrix $T$
(i.e.~with integer entries and determinant $\pm 1$),
having $\omega$ as an eigenvector with the associated eigenvalue
\[
  \lambda=\lambda(\Omega)>1
\]
(see Proposition~\ref{prop:matrixT} for an explicit construction).
Then, the iteration of the matrix $(T^{-1})\tp$
from an initial (\emph{``primitive''}) vector allows us
to generate any resonant sequence (see the definition~(\ref{eq:sqn})).

Next, we establish the \emph{main result} of this work,
providing two types of exponentially small asymptotic estimates
for the splitting, as $\eps\to0$,
related to the maximal distance of splitting and its transversality.
The first one is given by the maximum of $\abs{\M(\theta)}$, $\theta\in\T^2$.
On the other hand, we show that for most values of $\eps$
the function $\M(\theta)$ has simple zeros,
which correspond to transverse homoclinic orbits,
and for each zero $\theta^*$ we give an estimate for the
minimum eigenvalue (in modulus) of the matrix $\Df\M(\theta^*)$,
as a measure for the transversality of the splitting.
This generalizes the results of \cite{DelshamsG03,DelshamsG04}.

\begin{figure}[b!]
  \centering
  \subfigure{\label{fig:omega3a}
    \includegraphics[width=0.45\textwidth]{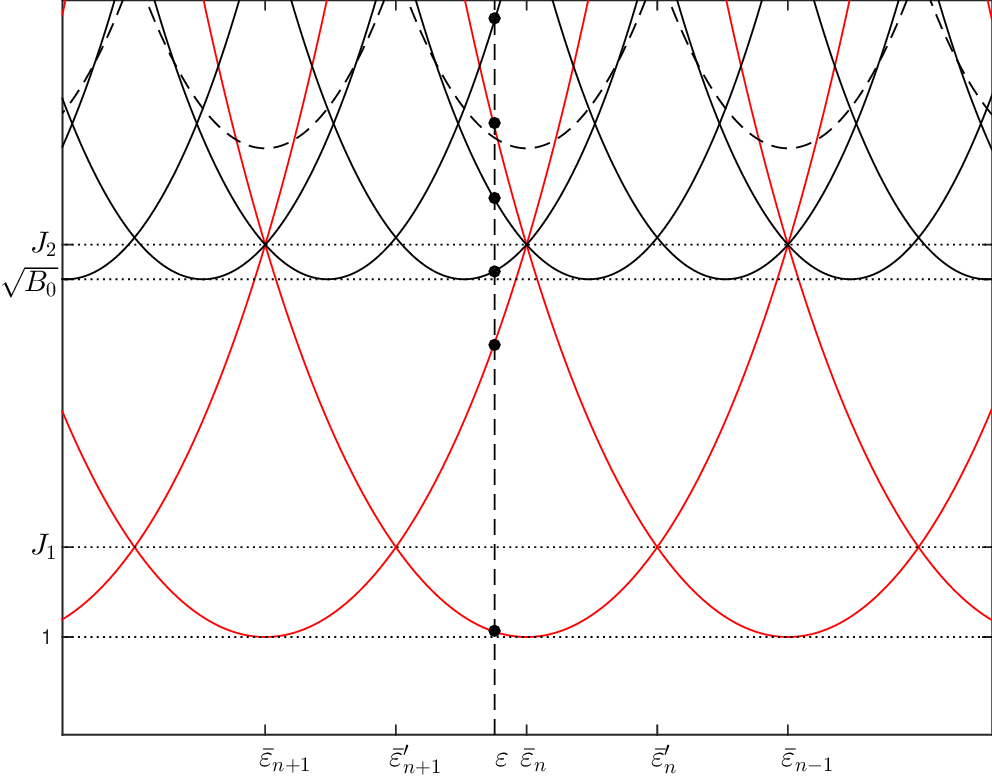}}
  \qquad
  \subfigure{\label{fig:omega3b}
    \includegraphics[width=0.48\textwidth]{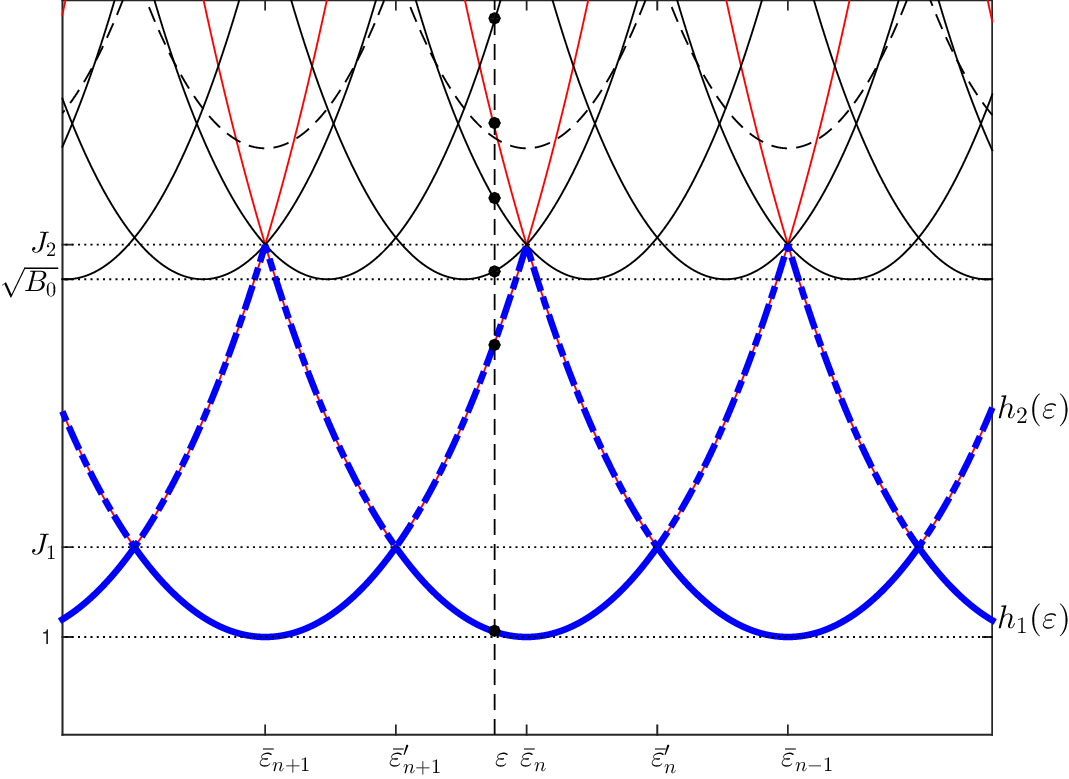}}
  \caption{\small\emph{
    Dependence on $\eps$ of the functions in the exponents,
    for the metallic ratio $\Omega=\ppcf{3}$ (the bronze ratio),
    using a~logarithmic scale for $\eps$,}
    \protect\\
    \subref{fig:omega3a}\emph{
      graphs of the functions $g^*_{s(q,n)}(\eps)$, associated to
      essential (the solid ones) and non-essential (the dashed ones)
      resonances~$s(q,n)$,
      the red ones correspond to the primary functions $\bg_n(\eps)$
      (see Section~\ref{sect:gn});}
    \protect\\
    \subref{fig:omega3b}\emph{
      graphs of the functions $h_1(\eps)$ and $h_2(\eps)$.}
    }
\label{fig:omega3}
\end{figure}

\begin{figure}[b!]
  \centering
  \vspace{2pt}
  \subfigure{\label{fig:omega13}
    \includegraphics[width=0.47\textwidth]{fig2a-splqua}}
  \qquad
  \subfigure{\label{fig:omega23}
    \includegraphics[width=0.47\textwidth]{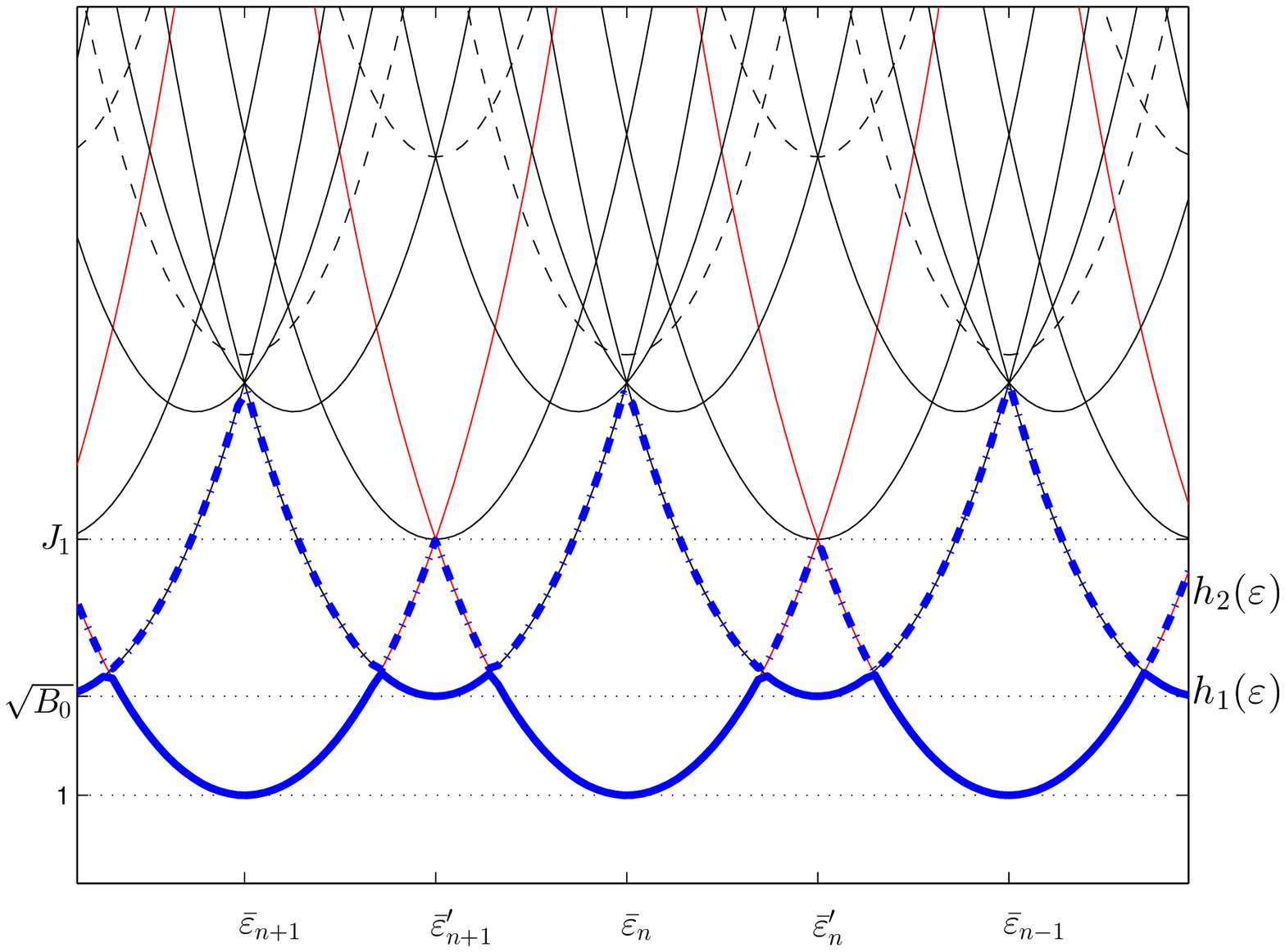}}
  \caption{\small\emph{
    Graphs of the functions $h_1(\eps)$ and $h_2(\eps)$
    for two metallic-colored ratios,}
    \protect\\
    \subref{fig:omega13}\emph{
      $\Omega=\ppcf{1,3}$ (a golden-colored ratio);}
    \qquad
    \subref{fig:omega23}\emph{
      $\Omega=\ppcf{2,3}$ (a silver-colored ratio).}
    }\label{fig:omega13-23}
\end{figure}

We stress that the dependence on $\eps$ of both asymptotic estimates
is given by the exponent 1/4, and by the functions $h_1(\eps)$ and $h_2(\eps)$,
which are periodic with respect to $\ln\eps$ and \emph{piecewise-smooth}
and, consequently, have a finite number of \emph{corners}
(i.e.~jump discontinuities of the derivative) in each period.
Some examples are shown in Figures~\ref{fig:omega3}--\ref{fig:omega13-23}
(where a \emph{logarithmic scale} for $\eps$ is used).
The oscillatory behavior of the functions $h_1(\eps)$ and $h_2(\eps)$
depends strongly on the arithmetic properties of $\Omega$
and, in fact, both functions can be explicitly constructed
from the continued fraction of~$\Omega$ (see Section~\ref{sect:h1h2}).
Below, in two additional results we establish more accurately
the behavior of the functions $h_1(\eps)$ and $h_2(\eps)$
in the simplest cases of 1-periodic and 2-periodic continued fractions.

For positive quantities, we use the notation $f\sim g$
if we can bound $c_1g\le f\le c_2g$ with constants
$c_1,c_2>0$ not~depending on $\eps$, $\mu$.

\begin{theorem}[\emph{main result}]\label{thm:main}
Assume the conditions described for the
Hamiltonian~\mbox{(\ref{eq:HamiltH}--\ref{eq:hf})},
with a quadratic frequency ratio $\Omega$,
that $\eps$ is small enough and that $\mu=\eps^r$, $r>3$.
Then, for the splitting function $\M(\theta)$ we have:
\btm
\item[\rm(a)]
$ \ds
\max_{\theta\in\T^2}\abs{\M(\theta)}
\sim \frac{\mu}{\eps^{1/2}}
\exp \pp{- \frac{C_0 h_1 (\eps)}{\eps^{1/4}}};
$
\item[\rm(b)]
the number of zeros $\theta^*$ of $\M(\theta)$ is $4\kappa$
with $\kappa(\eps)\ge1$ integer, and they are all simple,
for any $\eps$ except for a small neighborhood
of some transition values $\weps$,
belonging to a finite number of geometric sequences;
\item[\rm(c)]
at each zero $\theta^*$ of $\M(\theta)$,
the minimal eigenvalue of $\Df\M(\theta^*)$ satisfies
\[
  \abs{m^*} \sim \mu \eps^{1/4}
  \exp \pp{- \frac{C_0 h_2 (\eps)}{\eps^{1/4}}}.
\]
\etm
The functions $h_1(\eps)$ and $h_2(\eps)$, defined in~(\ref{eq:h1h2}),
are piecewise-smooth and $4\ln\lambda$-periodic in $\ln\eps$,
with $\lambda=\lambda(\Omega)$ as given in Proposition~\ref{prop:matrixT}.
In each period, the function $h_1(\eps)$ has at least 1~corner
and $h_2(\eps)$ has at least 2~corners.
They satisfy for $\eps>0$ the following bounds:
\[
  \min h_1(\eps)=1,
  \qquad
  \max h_1(\eps)\le J_1,
  \qquad
  \max h_2(\eps)\le J_2,
  \qquad
  h_1(\eps)\le h_2(\eps),
\]
with the constants
\beq\label{eq:defJ}
  J_1=J_1(\Omega):=\frac{1}{2}\p{\sqrt\lambda+\frac1{\sqrt\lambda}},
  \qquad
  J_2=J_2(\Omega):=\frac{1}{2}\p{\lambda+\frac1\lambda}.
\eeq
The corners of $h_1(\eps)$ are exactly the points $\ceps$ such that
$h_1(\ceps)=h_2(\ceps)$.
The corners of $h_2(\eps)$ are the same points~$\ceps$,
and the points $\weps$ where the results of~(b--c) do not apply.
The (integer) function $\kappa(\eps)$ is piecewise-constant
and $4\ln\lambda$-periodic in $\ln\eps$,
eventually with discontinuities at the transition points $\weps$.
On the other hand, $C_0=C_0(\Omega,\rho)$ is a positive constant
defined in~(\ref{eq:beta_gk}).
\end{theorem}

\bremarks
\item\label{rk:hf2}
If the function $h(x)$ in~(\ref{eq:hf}) is replaced by~$h(x)=\cos x -1$,
then the results of Theorem~\ref{thm:main} are valid for $\mu=\eps^r$
with $r>2$ (instead of $r>3$).
The details of this improvement are not given here,
since they work exactly as in~\cite{DelshamsG04}.
\item
As a consequence of the theorem, replacing $h_1(\eps)$ by
its upper bound $J_1>0$, we get the following \emph{lower bound} for
the maximal splitting distance:
\[
  \max_{\theta\in\T^2} |\M(\theta)|
  \ge \frac{c\mu}{\sqrt{\eps}}
  \exp \pp{- \frac{C_0J_1}{\eps^{1/4}}},
\]
where $c$ is a constant.
This may be enough, if our aim is only to prove the existence of splitting
of separatrices, without giving an accurate description for it.
\item
The results of Theorem~\ref{thm:main} can be partially generalized
if the frequency ratio $\Omega$ is a non-quadratic
\emph{number of constant type},
i.e.~whose continued fraction has bounded entries, but it is not periodic.
The numbers of constant type are exactly the ones such that
$\omega=(1,\Omega)$ satisfies
a Diophantine condition with the minimal exponent,
as in~(\ref{eq:DiophCond}).
This case has been considered in \cite{DelshamsGG14b},
where a function analogous to $h_1(\eps)$, providing an asymptotic estimate
for the maximal splitting distance, is defined. In this case,
this function is bounded but it is no longer periodic in $\ln\eps$.
\eremarks

For the simplest cases of continued fractions,
we can obtain more accurate information on the functions $h_1(\eps)$
and $h_2(\eps)$. As we show in Section~\ref{sect:contfract},
we can restrict ourselves to the case of
a \emph{purely periodic} continued fraction,
$\Omega=\ppcf{a_1,\ldots,a_m}$
(we assume that $0<\Omega<1$,
see Section~\ref{sect:contfract} for the notation).
In particular, we consider the following two cases:
\btm
\item the \emph{metallic ratios}:
\ $\Omega=\ppcf a=\dfrac{\dsqrt{a^2+4}-a}2$\,, \ with $a\ge1$;
\item the \emph{metallic-colored ratios}:
\ $\Omega=\ppcf{a,b}=\dfrac{\dsqrt{a^2b^2+4ab}-ab}{2a}$\,, \ with $1\le a<b$
\etm
(for the metallic-colored ratios, it is not necessary to consider
the case $a>b$, since $\ppcf{a,b}=\epcf a{b,a}$).

The metallic ratios,
which are limits of the sequence of quocients of consecutive terms
of generalized Fibonacci sequences, have often been considered
(see for instance \cite{Spinadel99,FalconP07}).
As some particular cases, we mention the \emph{golden}, \emph{silver}
and \emph{bronze} ratios:
\ $\ppcf{1}$, \ $\ppcf{2}$ \ and \ $\ppcf{3}$, \ respectively.
The metallic-colored ratios can be subdivided in several classes, such as:
\btm
\item the \emph{golden-colored ratios}: $\Omega=\ppcf{1,b}$, with $b\ge2$;
\item the \emph{silver-colored ratios}: $\Omega=\ppcf{2,b}$, with $b\ge3$;
\item the \emph{bronze-colored ratios}, etc.
\etm

For such types of ratios, we next provide descriptions of the
functions $h_1(\eps)$ and $h_2(\eps)$.
Additionally, in part~(b) we include a statement concerning
the exact number of critical points of the splitting function $\M(\theta)$,
for the case of metallic ratios.
Such results come from an accurate analysis of
the first and second essential dominant harmonics of the Melnikov potential,
studying whether they are both given by primary resonances for any $\eps$,
or they can be given by secondary resonances for some intervals of $\eps$.
We point out that a rigorous analysis of the role of the secondary resonances
becomes too cumbersome in some cases.
For this reason, although we give \emph{rigorous proofs}
for some results (the statements of Theorem~\ref{thm:corners}),
we provide \emph{numerical evidence} for other ones
(given below as ``Numerical Result~\ref{numres:corners2}'')
after checking them with intensive computations carried out
for a large number of frequency ratios
(see the proofs and justifications
in Section~\ref{sect:h1h2_metallic}).

\begin{theorem}\label{thm:corners}
Under the conditions of Theorem~\ref{thm:main}, we have:
\btm
\item[\rm(a)]
If the frequency ratio $\Omega$ is metallic,
the function $h_1(\eps)$ has exactly 1~corner $\ceps$ in each period,
satisfies \ $\max h_1(\eps)=h_1(\ceps)=J_1$,
\ and the distance between consecutive corners is exactly $4\ln\lambda$.
\item[\rm(b)]
If the frequency ratio $\Omega$ is golden-colored,
the function $h_2(\eps)$ has at least 3~corners in each period,
and satisfies \ $\min h_2(\eps)=J_1$, \ $\max h_2(\eps)<J_2$.
\item[\rm(c)]
If the frequency ratio $\Omega$ is metallic-colored but not golden-colored,
the function $h_1(\eps)$ has at least 2~corners in each period,
and satisfies \ $\max h_1(\eps)<J_1$.
\etm
\end{theorem}

\begin{numresult}\label{numres:corners2}
Under the conditions of Theorem~\ref{thm:main}, we have:
\btm
\item[\rm(a)]
If the frequency ratio $\Omega$ is golden-colored,
the function $h_1(\eps)$ has exactly 1~corner $\ceps$ in each period,
satisfies \ $\max h_1(\eps)=h_1(\ceps)=J_1$,
\ and the distance between consecutive corners is exactly $4\ln\lambda$.
\footnote{The result of part~(a) has been checked numerically
  for golden-colored ratios $\Omega=\ppcf{1,b}$, $2\le b\le10^6$.}
\item[\rm(b)]
If the frequency ratio $\Omega$ is metallic,
the function $h_2(\eps)$ has exactly 2~corners $\ceps$, $\weps$
in each period, satisfies \ $\min h_2(\eps)=h_2(\ceps)=J_1$,
\ $\max h_2(\eps)=h_2(\weps)=J_2$,
\ and the distance between consecutive corners is exactly~$2\ln\lambda$.
\ Moreover, the number of zeros $\theta^*$ of $\M(\theta)$ is exactly~4,
for any $\eps$ except for a small neighborhood
of the transition values $\weps$.
\footnote{The results of part~(b) have been
  checked numerically for metallic ratios $\Omega=\ppcf{a}$, $1\le a\le10^4$.}
\etm
\end{numresult}

As said before, the functions $h_1(\eps)$ and $h_2(\eps)$ can be
defined explicitly for any quadratic ratio $\Omega$,
from its continued fraction (see Section~\ref{sect:asympt_est}).
Such functions have piecewise expressions, which are simple in the case
of a metallic ratio, but in general they can be very complicated,
depending on the number of their corners in each period.
Nevertheless, we stress that numerical justifications
are required for results concerning infinite families of ratios,
such as the metallic or the golden-colored ones.
Instead, for particular frequency ratios the results
can be rigorously established.

\paragr{Organization of the paper}
We start in Section~\ref{sect:quadrfreq} with studying
the arithmetic properties of frequency vectors $\omega=(1,\Omega)$
with a quadratic ratio $\Omega$.
Such properties are closely related to the continued fraction of $\Omega$
(Section~\ref{sect:contfract}),
which allows us to construct the iteration matrices allowing us to
study the resonant properties of the vector $\omega$
(Sections~\ref{sect:resonantseq} and~\ref{sect:primary}),
and to provide accurate results for the cases of
metallic and metallic-colored ratios (Section~\ref{sect:metallic}),
mainly considered in this paper.
Next, in Section~\ref{sect:asympt_est} we find
an asymptotic estimate for the first and second dominant harmonics
of the splitting potential,
which allows us to define the functions $h_1(\eps)$ and $h_2(\eps)$
and study their general properties
(Sections~\ref{sect:gn} and~\ref{sect:h1h2}),
as well as the specific properties for
1-periodic and 2-periodic continued fractions
(Section~\ref{sect:h1h2_metallic}),
considered in Theorem~\ref{thm:corners}
and Numerical Result~\ref{numres:corners2}.
Finally, in Section~\ref{sect:technical} we provide rigorous bounds
of the remaining harmonics allowing us
to obtain asymptotic estimates for both the maximal splitting distance
and the transversality of the splitting,
as established in Theorem~\ref{thm:main}.

Finally, we introduce some notations that we use in this paper.
For positive quantities, we write $f\preceq g$ if we can bound $f\leq c\,g$
with some constant $c$ not depending on~$\eps$ and $\mu$.
In this way, we can write $f \sim g$ if $g\preceq f\preceq g$.
On the other hand, when comparing positive sequences $a_n$, $b_n$
we use an expression like ``\,$a_n\approx b_n$ as $n\to\infty$\,''
if \ $\ds\lim_{n\to\infty}(a_n/b_n)=1$,
\ and also ``\,$a_n\le b_n$ as $n\to\infty$\,''
if \ $\ds\limsup_{n\to\infty}(a_n/b_n)\le1$.

\section{Vectors with quadratic ratio}\label{sect:quadrfreq}

\subsection{Continued fractions of quadratic numbers}\label{sect:contfract}

It is well-known that any irrational number $0<\Omega<1$
has an infinite \emph{continued fraction}
\[
  \Omega=\npcf{a_1,a_2,a_3,\ldots}
  =\frac{1}{a_1+\dfrac{1}{a_2+\dfrac{1}{a_3+\cdots}}}\,,
  \qquad a_j\in\Z^+,\ j\ge1
\]
(notice that the integer part is $a_0=0$, hence
we have removed the entry `$0;$' from the notation).
Its entries $a_j$ are called the \emph{partial quotients}
of the continued fraction.
It is also well-known that the rational numbers
\ $\dfrac{p_j}{q_j}=[a_1,\ldots,a_j]$,\ $j\ge1$,
\ called the \emph{(principal) convergents} of $\Omega$,
provide successive best rational approximations to~$\Omega$.
Thus, if we consider the ``vector convergents'' $w(j):=(q_j,p_j)$,
we obtain approximations to
the direction of the vector $\omega=(1,\Omega)$
(see, for instance, \cite{Schmidt80} and \cite{Lang95}
as general references on continued fractions).

The convergents of a continued fraction are usually computed from the
standard recurrences
\beq\label{eq:convergents1}
  \begin{array}{lll}
    q_{-1}=0, &\ q_0=1,     &\qquad q_j=a_jq_{j-1}+q_{j-2},
  \\[4pt]
    p_{-1}=1, &\ p_0=a_0=0, &\qquad p_j=a_jp_{j-1}+p_{j-2},\quad j\ge1.
  \end{array}
\eeq
Alternatively, we can compute them in terms of products of unimodular matrices
\cite[Prop.~1]{DelshamsGG14b},
\beq\label{eq:convergents2}
  \mmatrix{q_j}{q_{j-1}}{p_j}{p_{j-1}}=A_1\cdots A_j,
  \qquad
  \mbox{where} \ A_i=\Tc(a_i):=\symmatrix{a_i}10.
\eeq
If we consider the first column, we can write
\ $w(j)=A_1\cdots A_jw(0)$.

An important tool in the study of continued fractions is the \emph{Gauss map}
$g:(0,1)\longrightarrow[0,1)$, defined as $g(x)=\pp{\dfrac1x}$,
where $\pp\cdot$ stands for the fractional part of any real number.
This map acts on a given continued fraction by removing the first entry:
for $\Omega=\npcf{a_1,a_2,a_3,\ldots}$,
we have $g(\Omega)=\npcf{a_2,a_3,\ldots}$.
We consider, for a given number $\Omega\in(0,1)$,
the sequence $(x_j)$ defined by
\beq\label{eq:gauss}
  x_0=\Omega
  \qquad
  x_j=g(x_{j-1}),\ j\ge1,
\eeq
which satisfies that $x_j\ne0$ for any $j$ if $\Omega$ is irrational.
It is clear that $x_j=\npcf{a_{j+1},a_{j+2},\ldots}$ for any $j$.

In our case of a \emph{quadratic irrational number} $\Omega$,
it is well-known that the continued fraction is \emph{eventually periodic},
i.e.~periodic starting at some partial quotient.
For an $m$-periodic continued fraction, we use the notation
\[
  \Omega=\epcf{b_1,\ldots,b_r}{a_1,\ldots,a_m}.
\]
In fact, as we see below we can restrict ourselves to the
numbers with \emph{purely periodic} continued fractions,
i.e.~periodic starting at the first partial quotient:
\ $\Omega=\ppcf{a_1,\ldots,a_m}$.
\ It is easy to relate such properties with the sequence $(x_j)$
defined by the Gauss map:
the continued fraction of $\Omega$ is eventually periodic
(hence, $\Omega$ is quadratic)
if and only if $x_{r+m}=x_r$ for some $r\ge0$, $m\ge1$,
and it is purely periodic if and only if $x_m=x_0$ for some $m\ge1$.

In the following proposition, which plays an essential role in
the results of this paper,
we see that for any given vector $\omega=(1,\Omega)$ 
with a quadratic ratio $\Omega$, there exists
a \emph{unimodular} matrix~$T=T(\Omega)$ having $\omega$
as an eigenvector with the associated eigenvalue $\lambda=\lambda(\Omega)>1$.
We show how we can construct both $T$ and $\lambda$,
directly from the continued fraction of $\Omega$.
Additionally, we show that applying the matrix $T$ to a convergent $w(j)$
we get the convergent~$w(j+m)$.

\begin{proposition}\label{prop:matrixT}
\ \btm
\item[\rm(a)]
Let $\Omega\in(0,1)$ be a quadratic irrational number
with a purely periodic continued fraction:
$\Omega=\ppcf{a_1,\ldots,a_m}$,
and consider the matrices $A_j=\Tc(a_j)$ as in~(\ref{eq:convergents2}).
Then, the matrix \ $T=A_1\cdots A_m$ \ is unimodular, and has
$\omega=(1,\Omega)$ as eigenvector with eigenvalue
$\lambda=\dfrac1{x_0x_1\cdots x_{m-1}}>1$,
where $(x_j)$ is the sequence defined by~(\ref{eq:gauss}).
Moreover, for the convergents $w(j)$ of $\Omega$
we have that \ $Tw(j)=w(j+m)$ \ for any $j\ge0$.
\item[\rm(b)]
Let $\wh\Omega$ be a quadratic irrational number
with a non-purely periodic continued fraction:
$\wh\Omega=\npcf{b_1,\ldots,b_r,\Omega}$, with~$\Omega$ as in~{\rm(a)},
and consider the matrices $B_j=\Tc(b_j)$, and \ $S=B_1\cdots B_r$.
\ Then, the matrix \ $\wh T=STS^{-1}$ \ is unimodular, and has
$\wh\omega=(1,\wh\Omega)$ as eigenvector with eigenvalue
$\lambda$ as in~{\rm(a)}.
Moreover, for the convergents $\hat w(j)$ of $\wh\Omega$
we have that \ $\wh T\hat w(j)=\hat w(j+m)$ \ for any $j\ge r$.
\etm
\end{proposition}

\proof
Using the construction of the sequence $(x_j)$ associated to $\Omega$,
we see that \ $\dfrac1{x_{j-1}}=a_j+x_j$, \ and we easily deduce the equality
\[
  \vect1{x_{j-1}}=x_{j-1}\,A_j\vect1{x_j},
  \quad n\ge1,
\]
Iterating this equality for $j=1,\ldots,m$ and using that $x_0=\Omega=x_m$,
we obtain
\[
  \vect1\Omega=x_0x_1\cdots x_{m-1}\;A_1A_2\cdots A_m\vect1\Omega,
\]
which proves that $T\omega=\lambda\omega$,
and it is clear that $T$ is unimodular.
To complete part~(a), using~(\ref{eq:convergents2}) and the
periodicity of the continued fraction we have
\[
  Tw(j)=A_1\cdots A_mA_1\cdots A_jw(0)=A_1\cdots A_{j+m}w(0)=w(j+m).
\]

With similar arguments we prove part~(b).
Indeed, using the sequence $(\hat x_j)$
associated to $\wh\Omega$, we see that
\[
  \vect1{\wh\Omega}
  =\hat x_0\hat x_1\cdots\hat x_{r-1}\;B_1B_2\cdots B_r\vect1\Omega,
\]
which says that the matrix $S$ provides a unimodular linear change
between the directions of the vectors $\omega$ and $\wh\omega$.
We deduce that $\wh T\hat\omega=\lambda\hat\omega$.
On the other hand, the matrix $S$ also provides a relation
between their respective convergents. Indeed, using~(\ref{eq:convergents2})
we see that, for $j\ge r$,
\[
  \hat w(j)
  =B_1\cdots B_rA_1\cdots A_{j-r}\hat w(0)
  =Sw(j-r)
\]
(notice that $\hat w(0)=w(0)=(1,0)$).
Then, using~(a) we deduce that
\[
  \wh T\hat w(j)=STw(j-r)=Sw(j)=\hat w(j+r).
\]
\qed

\bremarks
\item
In what concerns the contents of this paper, it is enough to consider
quadratic numbers with purely periodic continued fractions.
As we see from the proof of this proposition,
writing $S=\mmatrix{s_1}{s_2}{s_3}{s_4}$
we have the equality
\[
  \wh\Omega=\frac{s_3+s_4\Omega}{s_1+s_2\Omega}
  \qquad
  \mbox{with} \;\;\; s_1s_4-s_2s_3=\pm1,
\]
expressing the equivalence of the number $\wh\Omega$,
with an eventually periodic continued fraction,
with the number~$\Omega$ with a purely periodic one.
Then, it can be shown that our
main result (Theorem~\ref{thm:main})
applies to both numbers $\Omega$ and $\wh\Omega$ for $\eps$ small enough,
and we only need to consider the purely periodic case.
For instance, the results for the golden number $\Omega=\ppcf{1}$
also apply to the \emph{noble} numbers
$\wh\Omega=\epcf{b_1,\ldots,b_r}1$.
We point out that the treshold in $\eps$ of validity of the results,
not considered in this paper, would depend on the non-periodic part
of the continued fraction.
\item\label{rk:koch}
This proposition provides a particular case of
an algebraic result by Koch \cite{Koch99}, 
which also applies to higher dimensions:
for any given vector $\omega\in\R^n$ whose components
generate an algebraic number field of degree $n$,
there exists a unimodular matrix $T$ having $\omega$
as an eigenvector with the associated eigenvalue $\lambda$
of modulus greater than~1.
This result is usually applied in the context of renormalization theory,
since the iteration of the matrix~$T$ provides successive rational
approximations to the direction of the vector $\omega$
(see for instance \cite{Koch99,Lopesd02a}).
\eremarks

\subsection{Resonant sequences}\label{sect:resonantseq}

In this section and the next one,
we review briefly the technique developed in \cite{DelshamsG03}
for classifying the quasi-resonances of a given frequency vector
$\omega=(1,\Omega)$ whose ratio $\Omega$ is quadratic,
and study their relation with the convergents of the continued fraction
of $\Omega$.
A vector $k\in\Z^2\setminus\pp0$ can be considered
a quasi-resonance if $\scprod k\omega$ is small in modulus.
To determine the dominant harmonics of the Melnikov potential,
we can restrict to quasi-resonant vectors,
since the effect of vectors far enough from resonances can easily be bounded.

More precisely, we say that an integer vector $k\ne0$
is a \emph{quasi-resonance}
of $\omega$ if
\[
  \abs{\scprod k\omega}<\frac12\,.
\]
It is clear that any quasi-resonance can be presented in the form
\[
  k^0(q):=(-p,q),
  \qquad
  \mbox{with} \quad p=p^0(q):=\rint(q\Omega)
\]
(we denote $\rint(x)$ the closest integer to $x$).
Hence, we have the \emph{small divisors} $\scprod{k^0(q)}\omega=q\Omega-p$.
We denote by $\A$ the set of quasi-resonances $k^0(q)$ with $q\ge1$
(which can be assumed with no loss of generality).
We also say that $k^0(q)$ is an \emph{essential quasi-resonance}
if it is not a multiple of another integer vector
(if $p\ne0$, this means that $\gcd(q,p)=1$),
and we denote by $\A_0$ the set of essential quasi-resonances.

As said in Section~\ref{sect:contfract},
the matrix $T=T(\Omega)$ given by Proposition~\ref{prop:matrixT}
(in both cases of purely or non-purely periodic continued fractions)
provides approximations to the direction of $\omega=(1,\Omega)$.
Instead of $T$, we are going to use another matrix
providing approximations to the
orthogonal line~$\langle\omega\rangle^\bot$,
i.e.~to the quasi-resonances of $\omega$.
Notice the following simple but important equality:
\beq\label{eq:equalityU}
\scprod{(T^{-1})\tp k}\omega=\scprod k{T^{-1}\omega}
= \frac{1}{\lambda}\scprod k\omega,
\eeq
with $\lambda=\lambda(\Omega)$ as given by Proposition~\ref{prop:matrixT}.
With this in mind, for a quadratic ratio with an (eventually) $m$-periodic
continued fraction, we define the matrix
\beq\label{eq:defU}
  U=U(\Omega):=\sigma\,(T^{-1})\tp,
  \qquad
  \mbox{where}
  \ \sigma:=\det T=(-1)^m
\eeq
(the sign $\sigma$, which is not relevant, is introduced in order to have
a simpler expression in~(\ref{eq:abcd})).
It is clear from~(\ref{eq:equalityU}) that
if $k\in \A$, then also $Uk\in\A$. We say that the
vector $k=k^0(q)=(-p,q)$ is \emph{primitive}
if \ $k\in\A$ \ but $U^{-1}k\notin\A$.
\ If so, we also say that $q$~is a \emph{primitive integer},
and denote $\Pc$ the set of primitive integers, with $q\ge1$.
We deduce from~(\ref{eq:equalityU}) that $k$ is primitive if and only if
the following \emph{fundamental property} is fulfilled:
\beq\label{eq:primitive}
  \frac{1}{2\lambda} < |\langle k, \omega\rangle| < \frac{1}{2}\,.
\eeq
If a primitive $k^0(q)=(-p,q)$ is essential,
we also say that $q$ is an \emph{essential primitive integer},
and we denote $\Pc_0\subset\Pc$ the set of essential primitive integers.

Now we define, for each primitive vector $k^0(q)$,
the following \emph{resonant sequences} of integer vectors:
\beq\label{eq:sqn}
  s(q,n) := U^n k^0(q),
  \qquad n\ge0.
\eeq
It turns out that such resonant sequences cover the whole set of
vectors in $\A$, providing a classification for them.

\bremark
A resonant sequence $s(q,n)$ generated by an essential primitive $k^0(q)$
cannot be a multiple of another resonant sequence.
Indeed, in this case we would have $k^0(q)=c\,s(\tl q,n_0)$
with $c>1$ and $n_0\ge0$, and hence $k^0(q)$ would not be essential.
\eremark

Let us establish a relation between the resonant sequences $s(q,n)$, and
the convergents of $\Omega$. Alternatively to the convergents $w(j)=(q_j,p_j)$
considered in Section~\ref{sect:contfract},
we rather consider the \emph{``resonant convergents''}
(see also~\cite{DelshamsGG14b}),
\[
  v(j):=(-p_j,q_j).
\]
The next lemma shows that the action of the matrix $U$
defined in~(\ref{eq:defU}), on the vectors $v(j)$,
is analogous to the action of $T$ on the vectors $w(j)$
(which has been described in Proposition~\ref{prop:matrixT}).
This implies that the sequence of resonant convergents is divided
into $m$~of the resonant sequences defined in~(\ref{eq:sqn}).
We also see that the primitive vectors generating
such sequences are the $m$~first resonant convergents (belonging to $\A$).

\smallskip

\begin{lemma}\label{lm:jump}
\ \btm
\item[\rm(a)]
Let $\Omega$ be a quadratic number with
an (eventually) periodic continued fraction,
\ $\Omega=\epcf{b_1,\ldots,b_r}{a_1,\ldots,a_m}$ \ (with~$r\ge0$).
\ Then, we have
\[
  Uv(j)=v(j+m),
  \qquad j\ge r,
\]
and hence the sequence
of resonant convergents $v(j)$, for $j\ge r$,
is divided into $m$~resonant sequences.
\item[\rm(b)]
If $\Omega$ has a purely periodic continued fraction,
\ $\Omega=\ppcf{a_1,\ldots,a_m}$,
\ the primitive vectors among the resonant convergents are
\[
  \begin{array}{ll}
    v(1),\ldots,v(m)   &\quad\mbox{if} \ \ a_1=1;
  \\[4pt]
    v(0),\ldots,v(m-1) &\quad\mbox{if} \ \ a_1\ge2.
  \end{array}
\]
\etm
\end{lemma}

\proof
We use the following simple relation between the entries
of the matrices $T$ and $U$ (valid in both cases $r=0$ or~$r\ge1$):
\beq\label{eq:abcd}
  \mbox{if} \ \ T=\mmatrix abcd,
  \quad
  \mbox{then} \ \ U=\mmatrix d{-c}{-b}a,
\eeq
where we have taken into account that \ $\det T=(-1)^m$.
Then, the equality $Tw(j)=w(j+m)$, which holds for $j\ge r$,
is exactly the same as $Uv(j)=v(j+m)$, as stated in~(a),
by the relation between the vectors $v(j)$ and $w(j)$.
We have, as an immediate consequence, that the sequence
of resonant convergents $v(j)$ (for $j\ge r$)
is divided into $m$~resonant sequences.

To prove~(b), we first see that the small divisors associated
to the resonant convergents $v(j)$ satisfy the equality
\[
  q_j\Omega-p_j=(-1)^jx_0\cdots x_j,
  \qquad j\ge0,
\]
where $(x_j)$ is the sequence introduced in~(\ref{eq:gauss}).
This can easily be checked by induction, using the
recurrence~(\ref{eq:convergents1})
and the equality \ $\dfrac1{x_{j-1}}-a_j=x_j$.

If the continued fraction of $\Omega$ is purely periodic,
recalling the expression for $\lambda$ given
in Proposition~\ref{prop:matrixT}(a) and the
fundamental property~(\ref{eq:primitive}),
it is clear that a resonant convergent $v(j)$ is primitive if and only if
the following inequalities hold:
\[
  \frac{x_0\cdots x_{m-1}}2<x_0\cdots x_j<\frac12\,.
\]
Recall that $x_j\in(0,1)$ for any $j$. Using~(a), we see that such inequalities
can only be fulfilled by, at most, $m$~consecutive values of $j$.
For $a_1\ge2$, the first one is $j=0$ since \ $x_0=\Omega<1/2$,
\ and the last one is clearly $j=m-1$.
Instead, for $a_1=1$ the first one is $j=1$
since \ $x_0=\Omega>1/2$ \ and \ $x_0x_1=1-x_0<1/2$,
and the last one is $j=m$ since \ $x_m=x_0>1/2$.
\qed

\bremarks
\item\label{rmk:notriang}
The matrices $T$ and $U$ cannot be triangular,
i.e.~we have $b\ne0$ and $c\ne0$ in~(\ref{eq:abcd}).
Indeed, this would imply that the eigenvalue $\lambda$ is rational,
and hence the frequency ratio $\Omega$ would also be a rational number.
\item
The primitive resonant convergents given in part~(b) of this proposition
are all essential primitive vectors,
since all the convergents $p_j/q_j$ are reduced fractions
(as a consequence of the fact that the matrices in~(\ref{eq:convergents2})
are unimodular).
\eremarks

\subsection{Primary and secondary resonances}\label{sect:primary}

Now, our aim is to study which integer vectors $k$ fit best
the Diophantine condition~(\ref{eq:DiophCond}).
As in \cite{DelshamsG03}, we define the~\emph{``numerators''}
\beq\label{eq:numerators}
  \gamma_k := |\langle k, \omega\rangle|\cdot\abs k,
  \qquad
  k\in\Z^2\setminus\pp0
\eeq
where we use the norm $\abs\cdot=\abs\cdot_1$
(i.e.~the sum of absolute values).
As said in Section~\ref{sect:resonantseq}, we can restrict ourselves
to vectors $k=k^0(q)\in\A$ (with $q\ge1$),
and such vectors will be called primary or secondary resonances
depending on the size of $\gamma_k$.
We are also interested in studying the ``separation'' between
both types of resonances.

Recall that the matrix $T$ given by Proposition~\ref{prop:matrixT}
has $\omega=(1,\Omega)$ as an eigenvector with eigenvalue $\lambda>1$.
We consider a basis $\omega$, $v_2$ of eigenvectors of $T$,
where the second vector $v_2$ has the eigenvalue $\sigma/\lambda$
(of modulus~$<1$; recall that $\sigma=\det T$).
For the matrix $U$ defined in~(\ref{eq:defU}),
let $u_1$, $u_2$ be a basis of eigenvectors
with eigenvalues $\sigma/\lambda$ and~$\lambda$, respectively.
Writing the entries of the matrices $T$ and $U$ as in~(\ref{eq:abcd}),
it is not hard to obtain expressions for such eigenvalues and eigenvectors:
\bea
  \nonumber
  &\lambda=a+b\Omega,
  \qquad
  \dfrac\sigma\lambda=d-b\Omega,
\\[4pt]
  \label{eq:TU2}
  &v_2=(-b\Omega,c),
  \quad
  u_1=(c,b\Omega),
  \quad
  u_2=(-\Omega,1).
\eea
We also get the quadratic equations for the frequency ratio $\Omega$
and the eigenvalue $\lambda$:
\beq\label{eq:TU3}
  b\Omega^2=c-(a-d)\Omega,
  \qquad
  \lambda^2=(a+d)\lambda-\sigma.
\eeq
For any primitive integer $q$, recalling that we write $k^0(q)=(-p,q)$,
we define the quantities
\beq\label{eq:TU4}
    r_q:=\scprod{k^0(q)}\omega=q\Omega-p,
    \qquad
    z_q:=\scprod{k^0(q)}{v_2}=cq+bp\Omega.
\eeq

\bremark
As a consequence of the fact that $\Omega$ is an irrational number,
one readily sees that, if $q\ne\ol q$,
then $r_q\ne r_{\ol q}$ and $z_q\ne z_{\ol q}$
(in the latter case, using also that $c\ne0$,
as seen in remark~\ref{rmk:notriang} after Lemma~\ref{lm:jump}).
\eremark

The following proposition, whose proof is given in \cite{DelshamsG03}
(see also \cite{DelshamsGG14a} for a comparison with the case
of 3-dimensional cubic frequencies),
shows that the resonant sequences $s(q,n)$ defined in~(\ref{eq:sqn})
have a limit behavior: the sizes of the vectors $s(q,n)$
exhibit a \emph{geometric growth}, and the numerators $\gamma_{s(q,n)}$
tend to a \emph{``limit numerator''} $\gamma^*_q$,
as~$n\to\infty$.

\begin{proposition}\label{prop:quadrfreq}
Let $\omega=(1,\Omega)$ be a frequency vector whose ratio
$\Omega\in(0,1)$ is quadratic, and consider
the vectors $v_2$ and~$u_2$ as in~(\ref{eq:TU2}).
For any primitive integer $q\in\Pc$ (see~(\ref{eq:primitive})),
consider the quantities $r_q$ and $z_q$ defined in~(\ref{eq:TU4}).
Then, the resonant sequence $s(q,n)$ defined in~(\ref{eq:sqn}) satisfies:
\btm
\item[\rm(a)]
$|s(q,n)|= K_q \lambda^n + \Ord(\lambda^{-n})$,
\quad
where
\ $\ds K_q:=\abs{\frac{z_q}{\scprod{u_2}{v_2}}\,u_2}$;
\item[\rm(b)]
$\gamma_{s(q,n)}=\gamma^*_q + \Ord(\lambda^{-2n})$,
\quad
where
\ $\ds\gamma^*_q:=\lim_{n\to\infty}\gamma_{s(q,n)}=\abs{r_q}K_q$.
\etm
\end{proposition}

Using~(\ref{eq:TU2}--\ref{eq:TU4}),
we get the following alternative expression for the limit numerators:
\beq\label{eq:defdelta}
  \gamma^*_q=\frac{\Omega(1+\Omega)}{\abs{c+b\Omega^2}}\,\abs{\delta_q},
  \qquad
  \mbox{where}\ \ \delta_q:=\frac{r_qz_q}\Omega=cq^2-(a-d)qp-bp^2,
  \quad
  p=p^0(q).
\eeq
It is clear that $\delta_q\ne0$ and it is an integer.
We can select the minimal of the values $\abs{\delta_q}$
and, consequently, of the limit numerators $\gamma^*_q$,
which is reached by some concrete primitive $\wh q$.
We define
\beq\label{eq:quad_gamj0}
  \delta^*:=\min_{q\in\Pc}\abs{\delta_q}=\abs{\delta_{\wh q}}\ge1,
  \qquad
  \gamma^*:=\min_{q\in\Pc}\gamma^*_q=\gamma^*_{\wh q}>0.
\eeq
It is easy to see, as a consequence, that
\ $\ds\liminf_{\abs k\to\infty}\gamma_k=\gamma^*>0$.
\ Hence, any vector with quadratic ratio satisfies the
Diophantine condition~(\ref{eq:DiophCond}),
and we can consider $\gamma^*$ as the \emph{asymptotic Diophantine constant}.

As we see, all limit numerators $\gamma^*_q$ are multiple of
a concrete positive number.
An important consequence of this fact is that it
allows us to establish a classification of the vectors in $\A$.
We define \emph{the primary resonances} as the integer vectors
belonging to the sequence $s_0(n):=s(\wh q, n)$,
and \emph{secondary resonances} the vectors
belonging to any of the remaining sequences $s(q,n)$, $q\ne\wh q$
(recall that $\wh q$ is the primitive giving the minimum
in~(\ref{eq:quad_gamj0})).
We also introduce \emph{normalized numerators} $\tl\gamma_k$ and their limits
$\tl\gamma^*_q$, $q\in\Pc$, after dividing by $\gamma^*$,
and in this way $\tl\gamma^*_{\wh q}=1$.
We also define a value $B_0=B_0(\Omega)$ measuring the \emph{``separation''}
between the primary and the \emph{essential} secondary resonances:
\beq\label{eq:defB0}
  \tl\gamma_k:=\frac{\gamma_k}{\gamma^*}\,,
  \qquad
  \tl\gamma^*_q:=\frac{\gamma^*_q}{\gamma^*}=\frac{\abs{\delta_q}}{\delta^*}\,,
  \qquad
  B_0:=\min_{q\in\Pc_0\setminus\{\wh q\}}\tl\gamma^*_q.
\eeq

Using the fundamental property~(\ref{eq:primitive}) and the inequality
$\abs{p-q\Omega}<1/2$, we get the following lower bound for the
limit numerators,
which sligthly improves the analogous bound given in \cite{DelshamsG03}:
\beq\label{eq:lowerbound}
  \gamma^*_q>\frac{(1+\Omega)q-\alpha}{2 \lambda}\,,
  \qquad
  \alpha=\frac{\abs b\Omega(1+\Omega)}{2\abs{c+b\Omega^2}}\,.
\eeq

\bremarks
\item
Since the lower bound~(\ref{eq:lowerbound}) is increasing with respect to $q$,
it is enough to check a finite number of cases
in order to find the minimum in~(\ref{eq:quad_gamj0}).
\item\label{rmk:oneprimary}
We are implicitly assuming that the primitive integer $\wh q$ providing
the minimum in~(\ref{eq:quad_gamj0}) is unique.
In fact, we will show in Section~\ref{sect:metallic}
that this is true for the cases
of metallic or metallic-colored ratios $\Omega$
introduced in Section~\ref{sect:main}.
But in other cases, the minimum could be reached by two or more primitives
and, consequently, there could be two or more sequences of primary resonances.
For instance, it is not hard to check that for the ratio
$\Omega=\ppcf{1,2,2}$ there are two sequences of primary resonances.
\item\label{rmk:primary_ess}
Any primitive integer $\wh q$ generating a sequence of primary resonances
is essential. Indeed, if $\wh q$ is not essential, then we have
$k^0(\wh q)=c\,s(\ol q,n_0)$ with $c>1$ and $n_0\ge0$,
and therefore $s(\wh q,n)=c\,s(\ol q,n_0+n)$,
which implies by~(\ref{eq:numerators})
that $\gamma^*_{\wh q}=c^2\gamma^*_{\ol q}$,
and the minimum in~(\ref{eq:quad_gamj0})
would not be reached for $\wh q$.
\eremarks

Next, we show that the sequence of primary resonances is one (or more)
of the $m$~resonant sequences in which, by Lemma~\ref{lm:jump},
the resonant convergents are divided
if the continued fraction of $\Omega$ is $m$-periodic.
In fact, we can give a lower bound for the numerators of
all the remaining sequences.

\begin{lemma}\label{lm:primary}
For any primitive integer $q$ such that the vectors in the sequence $s(q,n)$
are not resonant convergents, its normalized numerator
satisfies $\tl\gamma^*_q>\sqrt5/2$.
\end{lemma}

\proof
We use some results in \cite[\S I.5]{Schmidt80}
(namely, Theorems~I.5B and I.5C),
concerning the properties of the convergents
of any irrational number.
On one hand, for an infinite number of convergents
the inequality $\abs{q_n\Omega-p_n}<1/(\sqrt5\,q_n)$ is satisfied; 
and on the other hand, if a given integer $q\ge1$ is not a convergent
and $p/q$ is a reduced fraction,
then $\abs{q\Omega-p}\ge1/2q$. 
To compare such results with our Diophantine condition~(\ref{eq:DiophCond}),
notice that $\abs{k^0(q)}=q+p\approx(1+\Omega)q$ as $q\to\infty$.

The first quoted result implies that, at least for one of
the resonant sequences $s(q,n)$ whose vectors are resonant convergents,
its limit numerator satisfies $\gamma^*_q\le(1+\Omega)/\sqrt5$.
By the second result, if a given resonant sequence $s(q,n)$
is generated by an essential primitive $q$
and its vectors are not resonant convergents,
then $\gamma^*_q\ge(1+\Omega)/2$.
This is also true if $q$ is not essential,
by the previous remark~\ref{rmk:primary_ess}.
Dividing the two bounds obtained, we get the lower bound $\sqrt5/2$
for the normalized limit $\tl\gamma^*_q$, when $q$ does not generate
a sequence of resonant convergents.
\qed

\subsection{Results for metallic and metallic-colored ratios}
\label{sect:metallic}

Now, we provide particular arithmetical results for the
(purely periodic) cases of a metallic ratio $\Omega=\ppcf{a}$,
and a metallic-colored ratio $\Omega=\ppcf{a,b}$,
introduced in Section~\ref{sect:main}.

\paragr{Metallic ratios}
Let us write, for a given $\Omega=\ppcf{a}$, $a\ge1$,
the matrix $T=T(\Omega)$ and the eigenvalue $\lambda=\lambda(\Omega)$,
as deduced from Proposition~\ref{prop:matrixT}(a),
and the matrix $U=U(\Omega)$ from~(\ref{eq:abcd}),
\beq\label{eq:Tmetallic}
  T=\symmatrix a10,
  \quad
  U=\symmatrix0{-1}a,
  \qquad
  \lambda=\dfrac1\Omega=a+\Omega.
\eeq
We also have from~(\ref{eq:TU3}) the quadratic equation
\beq\label{eq:Tmetallic2}
  \lambda^2=a\lambda+1.
\eeq
By Lemma~\ref{lm:jump}, all resonant convergents $v(j)$ belong to
a unique resonant sequence, whose primitive vector
is $v(1)=(-1,1)$ if $a=1$,
and $v(0)=(0,1)$ if $a\ge2$.
We deduce from Lemma~\ref{lm:primary} that this resonant sequence
provides the primary resonances: in both cases $\wh q=1$
and hence $s_0(n)=s(1,n)$.
In the next result, we compute the separation $B_0$,
defined in~(\ref{eq:defB0}), for \emph{all} metallic ratios,
providing in this way a sharp lower bound
for the normalized numerators of \emph{all} the essential secondary resonances.

\begin{proposition}\label{prop:metallic}
Let $\Omega=\ppcf{a}$, $a\ge1$, be a metallic ratio.
Then, the sequence of primary resonances is generated
by the primitive integer $\wh q=1$, and we have:
\[
  B_0=\tl\gamma^*_{q_1}=\left\{
    \begin{array}{ll}
        5 &\mbox{if}\ a=1,
      \\[3pt]
        a &\mbox{if}\ a\ge2,
    \end{array}
  \right.
  \qquad
  \mbox{for}
  \ \ q_1=\left\{
    \begin{array}{ll}
        7     &\mbox{if}\ a=1,
      \\[3pt]
        3     &\mbox{if}\ a=2,
      \\[3pt]
        a\pm1 &\mbox{if}\ a\ge3.
    \end{array}
  \right.
\]
\end{proposition}

\proof
We use the expression~(\ref{eq:defdelta}), taking into account the entries
of the matrix $T$ given in~(\ref{eq:Tmetallic}).
For the primary resonances, we have $\delta_1=\delta^*=1$,
and hence $\gamma^*=\dfrac{\Omega(1+\Omega)}{1+\Omega^2}$\,.
Dividing~(\ref{eq:lowerbound}) by $\gamma^*$
and using that $\lambda=1/\Omega$
we get, for the normalized numerators $\tl\gamma^*_q=\abs{\delta_q}$,
the following lower bound:
\[
  \tl\gamma^*_q>\frac{1+\Omega^2}2\,q-\frac\Omega4\,.
\]
If $a=1$ (the golden ratio), one checks that the second essential primitive
is $(-4,7)$ with $\tl\gamma^*_7=5$, and $\tl\gamma^*_q>5$ for $q\ge8$.
For $a\ge2$, assuming that $q>2/\Omega$ we get $\tl\gamma^*_q>a$.
Otherwise, if $q<2/\Omega$, since $p<q\Omega-1/2$ we get $p=p^0(q)<3/2$,
i.e.~$p=0$ or $p=1$. The only essential primitive with $p=0$ is $(0,1)$,
which gives the primary resonances, and for $p=1$ we have
an ``interval'' of primitives $(-1,q)$,
with \ $\dfrac{a+1}2\le q\le\dfrac{3a}2$ \ and \ $q\ne a$
\ (we have applied~(\ref{eq:primitive}) together with the fact
that $a<1/\Omega<a+1$). For such primitives, applying~(\ref{eq:defdelta})
we obtain \ $\delta_q=q^2-aq-1$, a quadratic polynomial in~$q$,
which is a increasing function for $q\ge a/2$,
with $\delta_{a\pm1}=\pm a$. This change of sign indicates that
$\tl\gamma^*_q=\abs{\delta_q}$ is minimal for $q=a\pm1$.
This argument is valid for $a=2$ (the silver ratio),
but in this case we must exclude $q=a-1$,
which lies outside the interval considered.
\qed

\paragr{Metallic-colored ratios}
Now, we consider $\Omega=\ppcf{a,b}$, $1\le a<b$.
Recall that, for $a=1$, this is called a golden-colored ratio;
we see below that our results are somewhat different for this particular case.
We have
\[
  T=\mmatrix{ab+1}ab1,
  \quad
  U=\mmatrix1{-b}{-a}{ab+1},
  \qquad
  \lambda
  =\frac1{1-a\Omega}=ab+1+a\Omega
\]
and, from~(\ref{eq:TU3}), the quadratic equation
\beq\label{eq:Tmetcol2}
  \lambda^2=(ab+2)\lambda-1.
\eeq
Applying Lemma~\ref{lm:jump}, we see that the resonant convergents $v(j)$
are divided into 2~resonant sequences, whose respective primitive vectors are 
\beq\label{eq:conv_metcol}
  \begin{array}{lll}
    v(1)=(-1,1), &v(2)=(-b,b+1) &\quad\mbox{if} \ \ a=1;
  \\[4pt]
    v(0)=(0,1),  &v(1)=(-1,a)   &\quad\mbox{if} \ \ a\ge2.
  \end{array}
\eeq
By Lemma~\ref{lm:primary}, one of the 2~sequences of resonant convergents
provides the primary resonances: $s_0(n)=s(\wh q,n)$.
We call \emph{the main secondary resonances} the vectors
in the second sequence, which we denote as $s_1(n):=s(\ol q,n)$.
Next, we find in Proposition~\ref{prop:metcol}
the value of $\tl\gamma^*_{\ol q}$,
i.e.~the limit numerator of the main secondary resonances, and this gives
a rigorous upper bound for the separation $B_0$.
In Numerical Result~\ref{numres:metcol2}
we present numerical evidence of the exact equality $B_0=\tl\gamma^*_{\ol q}$,
after checking it for a large number of ratios.
We point out that, for a given concrete frequency ratio $\Omega$,
the separation $B_0(\Omega)$ can be rigorously determined
since, by the lower bound~(\ref{eq:lowerbound}),
it is enough to consider the limit numerators $\tl\gamma^*_q$
for a finite number of essential primitive integers~$q$.

We also find
in Numerical Result~\ref{numres:metcol2}
the value of the \emph{``second separation''},
i.e.~the minimal normalized limit numerator
among the essential resonant sequences whose vectors
are not resonant convergents:
\beq\label{eq:defB1}
  B_1:=\min_{q\in\Pc_0\setminus\{\wh q,\ol q\}}\tl\gamma^*_q.
\eeq

\begin{proposition}\label{prop:metcol}
Let $\Omega=\ppcf{a,b}$, $1\le a<b$, be a metallic-colored ratio.
Then, the sequences of primary resonances and main secondary resonances
are generated, respectively, by primitive integers $\wh q$, $\ol q$
given by
\[
  \begin{array}{lll}
    \wh q=1, &\ \ol q=b+1 &\quad\mbox{if}\ \ a=1,
  \\
    \wh q=a, &\ \ol q=1   &\quad\mbox{if}\ \ a\ge2.
  \end{array}
\]
In both cases,
the separation satisfies
\[
  1<B_0\le\tl\gamma^*_{\ol q}=\frac ba\,.
\]
\end{proposition}

\proof
We use~(\ref{eq:defdelta}) in order to determine which
primitives~(\ref{eq:conv_metcol}) generate the sequence of primary resonances.
For $a=1$, we obtain $\delta_1=-1$ and $\delta_{b+1}=b$.
For $a\ge2$, we obtain $\delta_1=b$ and $\delta_a=-a$.
In both cases, the minimum (in modulus) is $\delta^*=a$,
which is reached for $\wh q=1$ if $a=1$, and $\wh q=a$ if $a\ge2$.
Then, we have $\ol q=b+1$ if $a=1$, and $\ol q=1$ if $a\ge2$,
and we obtain $\tl\gamma^*_{\ol q}=\abs{\delta_{\ol q}/\delta_{\wh q}}=b/a$,
which provides the upper bound: $B_0\le b/a$.
We also have the inequality $B_0>1$ as a consequence
of Lemma~\ref{lm:primary}.
\qed

\begin{numresult}\label{numres:metcol2}
In the notations of Proposition~\ref{prop:metcol},
the separation and the second separation are exactly
\footnote{The values of $B_0$ and $B_1$ have been checked numerically
  for all golden-colored ratios with \ $1=a<b\le10^6$,
  \ and for all metallic-colored ratios with \ $2\le a<b\le10^3$.}
\beq\label{eq:B01_metcol}
  B_0=\tl\gamma^*_{\ol q}=\frac ba\,,
  \qquad
  B_1=\left\{
  \begin{array}{ll}
    b+4               &\mbox{if}\ a=1,
  \\[5pt]
    \dfrac{(a-1)b+a}a &\mbox{if}\ a\ge2,
  \end{array}
  \right.
\eeq
which satisfy $1<B_0<B_1$.
\end{numresult}

\justif
Numerically, we can compute $B_1$ by bounding from below the limit numerators
$\tl\gamma^*_q$ for all the essential primitives $q\ne\wh q,\ol q$
(in view of~(\ref{eq:lowerbound}), only a finite number of primitives $q$
have to be considered).
We have checked that they all satisfy $\tl\gamma^*_q>b/a$
(at least for all the frequency ratios we have explored),
and hence we get $B_0=b/a$, and $B_1>B_0$.
We also obtain an expression for $B_1$,
given in~(\ref{eq:B01_metcol}) separately for the cases $a=1$ and $a\ge2$.
\qed

\bremark
The numerical explorations allow us to determine accurately the
primitive integers $q_2$ such that $B_1=\tl\gamma^*_{q_2}$,
i.e.~giving the minimum in~(\ref{eq:defB1}):
\[
  q_2=\left\{
  \begin{array}{ll}
    2          &\mbox{if}\ a=1,b=2,
  \\[3pt]
    3,\,9        &\mbox{if}\ a=1,b=3,
  \\[3pt]
    4,\,11       &\mbox{if}\ a=1,b=4,
  \\[3pt]
    5,\,8,\,13     &\mbox{if}\ a=1,b=5,
  \\[3pt]
    b+3        &\mbox{if}\ a=1,b\ge6,
  \\[3pt]
    2b+3       &\mbox{if}\ a=2,
  \\[3pt]
    a-1,\ ab+a+1 &\mbox{if}\ a\ge3.
  \end{array}
  \right.
\]
In each case, the primitive integers $q_2$ generate
the ``third most resonant'' sequences among the non-convergent ones
(i.e.~after the 2~sequences of resonant convergents).
Again, we stress that it is possible to obtain this kind of results thanks to
the lower bound~(\ref{eq:lowerbound}), which allows us to
carry out a finite number of computations for any given ratio~$\Omega$.
\eremark

\section{Searching for the asymptotic estimates}\label{sect:asympt_est}

In order to provide asymptotic estimates
for the splitting, we start with the first order approximation, given by the
Poincar\'e--Melnikov method.
Although our main result (Theorem~\ref{thm:main}) is stated
in terms of the splitting function $\M(\theta)=\nabla\Lc(\theta)$,
it is more convenient for us to work with
the (scalar) splitting potential $\Lc(\theta)$, whose
first order approximation is given by the Melnikov potential $L(\theta)$.

In this section, we provide the constructive part of the proof,
which amounts to find, for every sufficiently small $\eps$,
the first and the second dominant harmonics of the Fourier expansion
of the Melnikov potential $L(\theta)$,
with exponentally small asymptotic estimates for their size,
given by functions $h_1(\eps)$ and $h_2(\eps)$ in the exponents.
As a direct consequence of the arithmetic properties of quadratic ratios,
such functions are periodic with respect to $\ln\eps$.
We also study, from such arithmetic properties, whether the dominant harmonics
are given by primary resonances. This allows us to provide a more complete
description of the functions $h_1(\eps)$ and $h_2(\eps)$
in some particular cases
(Theorem~\ref{thm:corners}
and Numerical Result~\ref{numres:corners2}).

The final step in the proof of our main result
is considered in Section~\ref{sect:technical}.
It requires to provide bounds for the sum
of the remaining terms of the Fourier expansion of $L(\theta)$,
ensuring that it can be approximated by its dominant harmonics.
Furthermore, to ensure that the Poincar\'e--Melnikov
method~(\ref{eq:melniapproxM})
predicts correctly the size of the splitting in the singular case $\mu=\eps^r$,
one has to extend the results to the Melnikov function $\M(\theta)$
by showing that the asymptotic estimates of the dominant harmonics
are large enough to overcome the harmonics of the error term
in~(\ref{eq:melniapproxM}).
This step is analogous to the one done in \cite{DelshamsG04}
for the case of the golden number~$\Omega=\ppcf{1}$
(using the upper bounds for the error term provided in \cite{DelshamsGS04}).

\subsection{Estimates of the harmonics of the splitting potential}
\label{sect:gn}

We plug our functions $f$ and $h$, defined in~(\ref{eq:hf}),
into the integral~(\ref{eq:L}) and get the Fourier expansion
of the Melnikov potential,
where the coefficients can be obtained using residues:
\[
  L(\theta) = \sum_{k\in\Zc\setminus\pp0}
  L_k\,\cos(\langle k, \theta\rangle -\sigma_k),
  \qquad
  L_k = \frac{2\pi |\langle k, \omega_\eps\rangle|
  \,\ee^{-\rho\abs k}}{\sinh |\frac{\pi}{2}
  \langle k, \omega_\eps\rangle|}\,.
\]
We point out that the phases $\sigma_k$ are the same as in~(\ref{eq:hf}).
Using~(\ref{eq:omega_eps}) and taking into account the
definition of the numerators $\gamma_k$ in~(\ref{eq:numerators}),
we can present each coefficient $L_k=L_k(\eps)$, $k\in\Zc\setminus\pp0$,
in the form
\begin{equation}
\label{eq:alphabeta}
  L_k = \alpha_k\,\ee^{- \beta_k},
  \qquad
  \alpha_k(\eps) \approx \frac{4 \pi\gamma_k}{\abs k\sqrt{\eps}}\,,
  \quad
  \beta_k(\eps) =\rho\abs k + \frac{\pi \gamma_k}{2\abs k\sqrt{\eps}}\,,
\end{equation}
where an exponentially small term has been neglected in the denominator
of $\alpha_k$.
The most relevant term in this expression is $\beta_k$, which gives
the exponential smallness in $\eps$ of each coefficient,
and we will show that $\alpha_k$ provides a polynomial factor.
This says that, for any given $\eps$, the smallest exponents $\beta_k(\eps)$
provide the largest (exponentially small) coefficients $L_k(\eps)$
and hence the dominant harmonics. We are going to study the dependence
on $\eps$ of this dominance.

We start with providing a more convenient expression
for the exponents $\beta_k(\eps)$,
which shows that the smallest ones are~$\Ord(\eps^{-1/4})$
(this is directly related to the exponents $1/4$ in Theorem~\ref{thm:main}).
We introduce for any given $X$, $Y$ the function
\beq\label{eq:defG}
  G(\eps;X,Y):=
  \frac{Y^{1/2}}2\pq{\p{\frac\eps X}^{1/4}+\p{\frac X\eps}^{1/4}},
\eeq
which has its minimum at $\eps=X$
with $G(X;X,Y)=Y^{1/2}$ as the minimum value.
Notice that each function $G(\cdot;X,Y)$ is determined
by the point $(X,Y^{1/2})$.
Now, we define
\[
  g_k (\eps):=G(\eps;\eps_k,\tl\gamma_k),
  \qquad
  \eps_k:=\frac{D_0\tl\gamma_k^{\,2}}{\abs k^4}\,,
  \qquad
  D_0:=\p{\frac{\pi\gamma^*}{2\rho}}^2,
\]
and the functions $g_k(\eps)$ have their minimum at $\eps=\eps_k$,
with the minimal values $g_k(\eps_k) = \tilde{\gamma}_k^{1/2}$.
Recall that the asymptotic Diophantine constant $\gamma^*=\gamma^*_{\wh q}$
and the normalized numerators $\tl\gamma_k=\gamma_k/\gamma^*$
were introduced in~(\ref{eq:quad_gamj0}--\ref{eq:defB0}).
We deduce from~(\ref{eq:alphabeta}) that
\beq\label{eq:beta_gk}
  \beta_k(\eps) = \frac{C_0}{\eps^{1/4}}\,g_k (\eps),
  \qquad  
  C_0:=(2\pi\rho\gamma^*)^{1/2},
\eeq
and hence the lower bound
\ $\beta_k(\eps)\geq \dfrac{C_0\tl\gamma_k^{1/2}}{\eps^{1/4}}$\,.

Since we are interested in obtaining \emph{asymptotic estimates}
for the splitting and its transversality, rather than lower bounds,
we need to determine for any given $\eps$
the first and the second essential dominant harmonics,
which can be found among the smallest values $g_k(\eps)$.
To this aim it is useful to consider, for a given frequency ratio $\Omega$,
the graphs of the functions $g_k(\eps)$
associated to \emph{essential} quasi-resonances $k\in\A_0$
(recall that the notion of ``essentiality'' has been introduced
at the beginning of Section~\ref{sect:resonantseq}).
As an illustration, such graphs are shown
in Figure~\ref{fig:omega3}\subref{fig:omega3a}
for a concrete example (the bronze ratio $\Omega=\ppcf{3}$),
using a logarithmic scale for $\eps$.
Other examples are shown in Figures~\ref{fig:omega13-23}--\ref{fig:omega122}.
The periodicity which can be noticed from the graphs can easily be explained
from the classification of the integer vectors into resonant sequences
(recall their definition in~(\ref{eq:sqn})).
Indeed, for $k=s(q,n)$ belonging to a concrete resonant sequence,
using the approximations for $|s(q,n)|$ and
$\gamma_{s(q,n)}$ given by Proposition~\ref{prop:quadrfreq},
we obtain the following approximations as~$n\to\infty$,
\beq\label{eq:gqn}
  g_{s(q,n)}(\eps) \approx g^*_{s(q,n)}(\eps)
  :=G(\eps;\eps^*_{s(q,n)},\tl\gamma^*_q),
  \qquad
  \eps_{s(q,n)} \approx \eps^*_{s(q,n)}
  :=\frac{D_0(\tl\gamma^*_q)^2}{K_q^{\,4}\lambda^{4n}}\,,
\eeq
which motivates the use of a logarithmic scale.
We point out that the graphs shown
in Figure~\ref{fig:omega3}\subref{fig:omega3a}
do not correspond to the true functions $g_{s(q,n)}(\eps)$,
but rather to the approximations $g^*_{s(q,n)}(\eps)$,
which satisfy the following scaling property:
\beq\label{eq:scaling}
  g^*_{s(q,n+1)}(\eps)=g^*_{s(q,n)}(\lambda^4\eps).
\eeq
This gives, for any resonant sequence, the mentioned periodicity:
the graph of $g^*_{s(q,n+1)}$ is a translation of $g^*_{s(q,n)}$,
to distance~$4\ln\lambda$.
For non-essential resonant sequences, whose vectors do not belong to $\A_0$,
we see that, if $s(q,n)=c\,s(\ol q,n_0+n)$ with $c>1$ and $n_0\ge0$,
then
\beq\label{eq:nonessential1}
  g^*_{s(q,n)}(\eps)=c\,g^*_{s(\ol q,n_0+n)}(\eps)
\eeq
(see also the remark~\ref{rmk:primary_ess} just before Lemma~\ref{lm:primary}).

In order to study the dependence on $\eps$ of the most dominant harmonics,
it is useful to study the intersections between the graphs of different
functions $g^*_k(\eps)$, since this gives the values of $\eps$ at which
a change in the dominance may take place.
In the next lemma, we consider the graphs associated to two different
quasi-resonances $k,\ol k\in\A$,
and we show that only two situations are posible:
they do not intersect (which says that one of them always dominates
the other one), or they intersect transversely in a unique point
(and in this case a unique change in the dominance takes place).

\begin{lemma}\label{lm:ZW}
Let $k,\ol k\in\A$, with $k\ne\ol k$,
given by $k=s(q,n)$ and $\ol k=s(\ol q,\ol n)$,
and assume that \ $\tl\gamma^*_q\le\tl\gamma^*_{\ol q}$\,.
\ Denoting \ $Z=\p{\eps^*_{\ol k}/\eps^*_k}^{1/4}$ \ and
\ $W=\p{\tl\gamma^*_{\ol q}/\tl\gamma^*_q}^{1/2}$,
\ the graphs of the functions $g^*_k(\eps)$ and $g^*_{\ol k}(\eps)$
intersect if and only if \ $Z<1/W$ \ or \ $Z>W$.
\ If so, the intersection is unique and transverse,
and takes place at \ $\eps=\eps^*_k\cdot\p{\dfrac{Z(WZ-1)}{Z-W}}^2$.
\end{lemma}

\proof
First of all, we show that $g^*_k$ and $g^*_{\ol k}$
cannot be the same function.
By the definition~(\ref{eq:defG}), if $g^*_k=g^*_{\ol k}$ then we have
\ $\tl\gamma^*_q=\tl\gamma^*_{\ol q}$
\ and \ $\eps^*_k=\eps^*_{\ol k}$\,.
\ The latter equality implies that
\ $K_q\lambda^n=K_{\ol q}\lambda^{\ol n}$.
\ Using the expressions given in Proposition~\ref{prop:quadrfreq},
we get the equalities
\ $\abs{r_qz_q}=\abs{r_{\ol q}z_{\ol q}}$
\ and \ $z_q\lambda^n=z_{\ol q}\lambda^{\ol n}$.
We deduce that \ $\abs{r_{\ol q}/r_q}=\lambda^{\ol n-n}$,
but we have \ $\abs{r_q},\abs{r_{\ol q}}\in(1/2\lambda,1/2)$
\ by the fundamental property~(\ref{eq:primitive}).
This says that $n=\ol n$ and hence $z_q=z_{\ol q}$.
As seen in the remark next to the definition~(\ref{eq:TU4}),
we also get $q=\ol q$, which contradicts the assumption~$k\ne\ol k$.

Now, introducing the variable $\zeta=\p{\eps/\eps^*_k}^{1/4}>0$, we define
\[
  f_1(\zeta):=\frac12\p{\zeta+\frac1\zeta}
  =\frac{g^*_k(\eps)}{\p{\tl\gamma^*_q}^{1/2}}\,,
  \qquad
  f_2(\zeta):=\frac W2\p{\frac\zeta Z+\frac Z\zeta}
  =\frac{g^*_{\ol k}(\eps)}{\p{\tl\gamma^*_q}^{1/2}}\,,
\]
with $W\ge1$ by hypothesis, and it is clear from the above analysis
that we cannot have $W=Z=1$.
It is straightforward to check that the graphs of $f_1$ and $f_2$
can intersect only once, transversely, at \ $\zeta^2=\dfrac{Z(WZ-1)}{Z-W}$\,.
\ Such an intersection occurs if and only if \ $Z<1/W$ \ or \ $Z>W$.
\ Then, we get the result after translating
from $\zeta$ to the original variable~$\eps$.
\qed

The sequence of primary resonances $s_0(n)=s(\wh q, n)$,
introduced in Section~\ref{sect:primary} plays an important role,
since they give the smallest minimum values among the functions $g^*_k(\eps)$,
and hence they will provide the most dominant harmonics,
at least for $\eps$ close to such minima.
With this fact in mind, and recalling that $\tl\gamma^*_{\wh q}=1$,
we denote
\bea
  \label{eq:gn}
  &&\bg_n(\eps)
  :=g^*_{s_0(n)}(\eps)
  =G(\eps;\beps_n,1)
  =\frac12\pq{\p{\frac\eps{\beps_n}}^{1/4}+\p{\frac{\beps_n}\eps}^{1/4}},
\\
  \nonumber
  &&\beps_n:=\eps^*_{s_0(n)}=\frac{D_0}{K_{\wh q}^{\,4}\lambda^{4n}}\,.
\eea

To study the periodicity with respect to $\ln\eps$,
we introduce intervals $\I_n$ whose ``length'' (in the logarithmic scale)
is~$4\ln\lambda$, centered at $\beps_n$,
and the left and right ``halves'' of such intervals,
\beq\label{eq:In}
  \I_n:=\pq{\peps_{n+1},\peps_n}=\I^+_n\cup\I^-_n,
  \qquad
  \I^+_n:=\pq{\peps_{n+1},\beps_n},
  \quad
  \I^-_n:=\pq{\beps_n,\peps_n},
\eeq
where \ $\ds\peps_n:=\sqrt{\beps_n\beps_{n-1}}=\lambda^2\beps_n$
\ are the geometric means of the sequence~$\beps_n$.
For a given $n\ge1$, it is easy to determine the behavior
of the functions~(\ref{eq:gn}):
for $\eps\in\I^+_n$, the value of the function $\bg_n(\eps)$
decreases from $J_1$ to 1,
the value of the function $\bg_{n+1}(\eps)$ increases from $J_1$ to $J_2$,
and we have $\bg_m(\eps)\ge J_2$ if $m\ne n,n+1$
(recall that the values $J_1$ and $J_2$ were defined in~(\ref{eq:defJ})).
A~symmetric result holds for $\eps\in\I^-_n$
with the functions $\bg_n(\eps)$ and $\bg_{n-1}(\eps)$
(see the red graphs in Figure~\ref{fig:omega3}\subref{fig:omega3a}
for an illustration).

\subsection{Dominant harmonics of the splitting potential}\label{sect:h1h2}

In this section, we introduce the functions $h_1(\eps)$ and $h_2(\eps)$
appearing in the exponents in Theorem~\ref{thm:main},
as the first and second minima, for any given $\eps$,
of the values $g^*_k(\eps)$ among the essential quasi-resonances $k\in\A_0$.
We study some of the properties of $h_1(\eps)$ and $h_2(\eps)$,
which hold for an arbitrary quadratic ratio $\Omega$.
In Section~\ref{sect:h1h2_metallic}, we put emphasis on the dependence
of such functions on the continued fraction of the frequency ratio $\Omega$,
giving a more accurate description of them,
for the cases of metallic and metallic-colored ratios,
whose arithmetic properties have been studied in Section~\ref{sect:metallic}.

\begin{figure}[b!]
  \centering
  \subfigure{\includegraphics[width=0.45\textwidth]{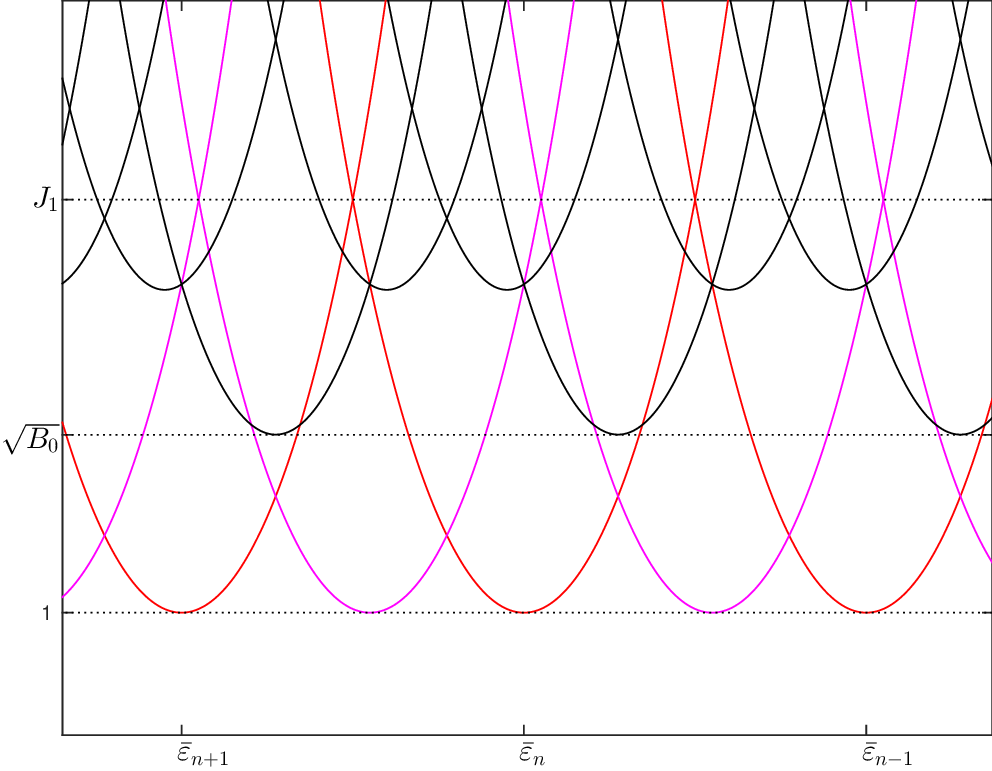}}
  \qquad
  \subfigure{\includegraphics[width=0.48\textwidth]{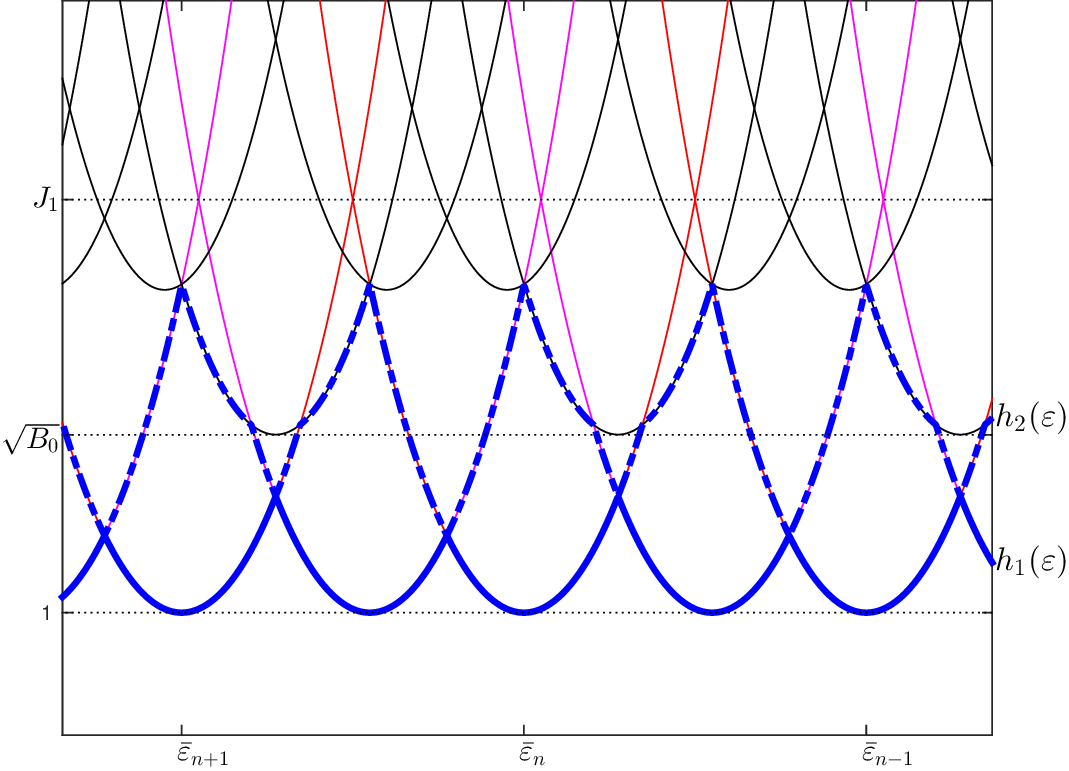}}
  \caption{\small\emph{
    Graphs of the functions $g^*_{s(q,n)}(\eps)$, $h_1(\eps)$ and $h_2(\eps)$
    for $\Omega=\ppcf{1,2,2}$:
    the two sequences of primary resonances correspond
    to the red and magenta graphs; the non-essential resonances
    are not represented.}}\label{fig:omega122}
\end{figure}

Namely, we provide information on the minimum and maximum values
of the functions $h_1(\eps)$ and $h_2(\eps)$, and show that they
are piecewise-smooth and $4\ln\lambda$-periodic in $\ln\eps$,
and give lower bounds for the number of their corners
(i.e.~jump discontinuities of the derivative)
in any period, say the interval \ $\I_n=[\peps_{n+1},\peps_n]$,
\ or rather the semi-open interval \ $(\peps_{n+1},\peps_n]$,
\ to avoid repetitions if the endpoints are corners.
In fact, such properties are clear from
Figures~\ref{fig:omega3}--\ref{fig:omega122},
for the concrete frequency rations considered there,
but we are going to show that they hold for
an arbitrary quadratic ratio $\Omega$.

Previously to this, let us define two functions analogous
to $h_1(\eps)$ and $h_2(\eps)$,
but taking into account \emph{only the primary} resonances:
\beq\label{eq:bh1bh2}
  \bh_1(\eps):=\min_n\bg_n(\eps)=g^*_{N_1}(\eps),
  \qquad
  \bh_2(\eps):=\min_{n\ne N_1}g^*_n(\eps)=g^*_{N_2}(\eps),
\eeq
with $N_i=N_i(\eps)$. In other words, the two dominant harmonics
among the primary resonances
correspond to
\[
  \bS_i=\bS_i(\eps)=s_0(N_i),
  \qquad
  i=1,2.
\]
On each concrete interval $\I_n$ (see the definition~(\ref{eq:In}))
one readily sees, from the properties
described in the last part of Section~\ref{sect:gn},
what primary resonances provide the first and second minima:
\ $N_1(\eps)=n$ for $\eps\in\I_n$,
\ and \ $N_1(\eps)=n\pm1$ for $\eps\in\I^\pm_n$.
\ It is also clear that the functions $\bh_1(\eps)$ and $\bh_2(\eps)$
are piecewise-smooth and $4\ln\lambda$-periodic.
In each period,
the function $h_1(\eps)$ has exactly 1~corner (at $\peps_n$),
and $h_2(\eps)$ has exactly 2~corners (at $\peps_n$ and $\beps_n$).
Moreover, we have
\bean
  &&\min\bh_1(\eps)=\bh_1(\beps_n)=1,
  \qquad
  \max\bh_2(\eps)=\bh_2(\beps_n)=J_2,
\\
  &&\max\bh_1(\eps)=\min\bh_2(\eps)=\bh_1(\peps_n)=\bh_2(\peps_n)=J_1
\eean
(see also Figure~\ref{fig:omega3}\subref{fig:omega3b}
for an illustration).

Now, we define the functions $h_i(\eps)$ as the 
minimal values of the functions $g^*_k(\eps)$
among \emph{all essential} quasi-resonances,
and we denote $S_i=S_i(\eps)$ the integer vectors $k$
at which such minima are reached:
\beq\label{eq:h1h2}
  \begin{array}{ll}
    \ds h_1(\eps):=\min_{k\in\A_0}g^*_k(\eps)=g^*_{S_1}(\eps),
    &\quad
    \ds h_2(\eps):=\min_{k\in\A_0\setminus\pp{S_1}}g^*_k(\eps)=g^*_{S_2}(\eps),
  \\[10pt]
    \ds h_3(\eps):=\min_{k\in\A_0\setminus\pp{S_1,S_2}}g^*_k(\eps)=g^*_{S_3}(\eps).
  \end{array}
\eeq
It is clear that $h_i(\eps)\le\bh_i(\eps)$ for any $\eps$ and $i=1,2$.
In order to provide an accurate description of the splitting and its
transversality, we have to study whether the equality between
the above functions can be established for any value of $\eps$,
or at least for some intervals of $\eps$.
This amounts to study whether the dominant harmonics
can be always found among the primary resonances ($S_i=\bS_i$)
or, on the contrary, secondary resonances have to be taken into account.

In fact, the properties described above for the functions $\bh_i(\eps)$
are partially generalized to the functions $h_i(\eps)$ in the next proposition,
which corresponds to some parts of the statement of Theorem~\ref{thm:main},
concerning such functions.
Recall that the values $J_1$ and $J_2$ were defined in~(\ref{eq:defJ}).

\begin{proposition}\label{prop:h1h2}
The functions $h_1(\eps)$ and $h_2(\eps)$
are piecewise-smooth, $4\ln\lambda$-periodic in $\ln\eps$.
In each period, the function $h_1(\eps)$ has at least 1~corner
and $h_2(\eps)$ has at least 2~corners.
They satisfy for $\eps>0$ the following bounds:
\[
  \min h_1(\eps)=1,
  \qquad
  \max h_1(\eps)\le J_1,
  \qquad
  \max h_2(\eps)\le J_2,
  \qquad
  h_1(\eps)\le h_2(\eps).
\]
The corners of $h_1(\eps)$ are exactly the points $\ceps$ such that
$h_1(\ceps)=h_2(\ceps)$.
The corners of $h_2(\eps)$ are the same points~$\ceps$,
and the points $\weps$ where $h_2(\weps)=h_3(\weps)$.
\end{proposition}

\proof
First of all, it is clear that the functions $h_1$ and $h_2$
are $4\ln\lambda$-periodic,
as we see from the scaling property~(\ref{eq:scaling}).
Then, we can restrict ourselves to a concrete interval, say $\I_1$.
Recalling also that
$h_i(\eps)\le\bh_i(\eps)\le J_i$ for any $\eps$ and $i=1,2$,
the minimum in the definition~(\ref{eq:h1h2}) of $h_i$
can be restricted to the integer vectors $k=s(q,n)\in\A_0$
such that the graph of the function $g^*_{s(q,n)}$
visits (the interior of) the rectangle \ $\I_1\times[1,J_i]$.
\ We are going to show that this is possible only
for a finite number of integer vectors.

Indeed, recalling~(\ref{eq:defG}), if the graph of a function $G(\eps;X,Y)$
visits \ $\I_1\times[1,J_1]$ \ then
\ $Y^{1/2}<J_1$ \ and \ $\beps_2<X<\beps_0$;
\ and if it visits \ $\I_1\times[1,J_2]$
\ then \ $Y^{1/2}<J_2$ \ and \ $\peps_3<X<\peps_0$.
\ For the function $g^*_{s(q,n)}$, defined in~(\ref{eq:gqn}),
we have to consider $Y=\tl\gamma^*_q$ and $X=\eps^*_{s(q,n)}$.
By the lower bound~(\ref{eq:lowerbound}), it is clear that only a finite
number of functions $g^*_{s(q,n)}$ can visit the rectangle
\ $\I_1\times[1,J_i]$, \ $i=1,2$.
\ This implies, by Lemma~\ref{lm:ZW}, that only a finite number of
(transverse) intersections between the graphs of $g^*_{s(q,n)}$
can take place inside the rectangles \ $\I_1\times[1,J_i]$.

We deduce from the above considerations that the functions
$h_1$ and $h_2$ are piecewise-smooth.
Indeed, we can consider a partition of $\I_1$ into subintervals such that,
for $\eps$ belonging to (the interior of) each subinterval,
the function~$h_1$ coincides
with only one of the functions $g^*_{s(q,n)}$,
i.e.~the dominant harmonic is given by $S_1(\eps)=s(q,n)$,
which remains constant on this subinterval.
At each endpoint of such subintervals, a change in the dominant harmonic
takes place, i.e.~$S_1$ has a jump discontinuity.
By Lemma~\ref{lm:ZW}, the endpoints of the subintervals correspond
to transverse intersections between
the graphs of different functions $g^*_{s(q,n)}$,
which give rise to corners $\ceps$ of $h_1$.
A similar argument applies to the function $h_2$,
with a different partition, associated to the changes
of the second dominant harmonic $S_2(\eps)$.
In fact, the values $\ceps$ are the points where $h_1(\ceps)=h_2(\ceps)$,
and they are corners of both functions $h_1$ and $h_2$.
In the same way, the function $h_2$ has additional corners $\weps$
at the points where $h_2(\weps)=h_3(\weps)$.

Finally, we provide a lower bound for the number of
corners $\ceps$, $\weps$ in a given period
(if the endpoints of a period are corners, we count them as one single corner).
Since $\bg_1(\beps_1)=1$, we have $S_1(\eps)=\bS_1(\eps)=s_0(1)$
in some neighborhood of $\beps_1\in\I_1$.
Analogously, we have $S_1(\eps)=\bS_1(\eps)=s_0(2)$
in some neighborhood of $\beps_2\in\I_2$,
which implies the existence of at least one corner of $h_1$
with $\beps_2<\ceps<\beps_1$ and, consequently, in any given period.
On the other hand, if $\ceps<\tl\eps$ are two consecutive corners
of $h_1$ (and $h_2$), there exists
at least one additional corner $\weps$ of $h_2$ (and $h_3$),
since $S_1(\eps)$ and $S_2(\eps)$
cannot be simultaneously constant in the interval $[\ceps,\tl\eps]$
(this would imply that $g^*_{S_1}$ and $g^*_{S_2}$
intersect at both points $\ceps$ and $\tl\eps$,
which is not possible by Lemma~\ref{lm:ZW}).
\qed

\bremarks
\item
We can also deduce from the proof of this proposition some useful
properties of the functions $S_i=S_i(\eps)$, giving the dominant harmonics.
Namely, each function $S_i(\eps)$ is ``piecewise-constant'',
with jump discontinuities at the corners of $h_i(\eps)$.
Moreover, the asymptotic behavior of the functions $S_i(\eps)$ as $\eps\to0$
turns out to be polynomial:
\beq\label{eq:estimS}
  \abs{S_i(\eps)} \sim \frac1{\eps^{1/4}}\,.
\eeq
Indeed, the first dominant harmonic belongs to some resonant sequence:
we can write $S_1(\eps)=s(q,N)$
for some $q=q(\eps)$, and for $N=N(\eps)$ such that the value $\eps^*_{s(q,N)}$
is the closest to $\eps$, among the sequence $\eps^*_{s(q,n)}$, $n\ge0$.
Recalling~(\ref{eq:gqn}) and the estimate $\abs{s(q,N)}\sim\lambda^N$ given
in Proposition~\ref{prop:quadrfreq}(a), we get~(\ref{eq:estimS}).
An analogous argument holds for $S_2(\eps)$, possibly replacing $N$ by $N\pm1$,
and possibly belonging to a different resonant sequence $s(q,\cdot)$.
Notice that it is not necessary to include $q$
in the estimate~(\ref{eq:estimS})
(in spite of the fact that $K_q$ and $\tl\gamma^*_q$
appear in the expression~(\ref{eq:gqn})),
since by the arguments in the above proof (Proposition~\ref{prop:h1h2})
only a finite number of resonant sequences $s(q,\cdot)$ can be involved.
\item
A more careful look at the arguments of the previous remark,
says that, if the dominant harmonics in a given interval~$\I_n$ are known,
then in view of the scaling property~(\ref{eq:scaling})
the dominant harmonics in the interval $\I_{n+1}$
are the next vectors in the respective resonant sequences:
\[
  S_i(\eps)=U\,S_i(\lambda^4\eps)
\]
(recall that the matrix $U$ appears in the definition of the resonant sequences
in~(\ref{eq:sqn})).
\item
Although we implicitly assume that there exists only one sequence of primary
resonances (see remark~\ref{rmk:oneprimary} before Lemma~\ref{lm:primary}),
it is not hard to adapt the definitions and results to the case of
two or more sequences of primary resonances.
In this case, we would choose one of such sequences
as ``the'' sequence $s_0(n)$, when the functions $\bg_n(\eps)$
are defined in~(\ref{eq:gn}).
As an example, we show in Figure~\ref{fig:omega122} the graphs of $h_1(\eps)$
and $h_2(\eps)$ for the ratio $\Omega=\ppcf{1,2,2}$, with
two sequences of primary resonances.
\eremarks

\subsection{Dominant harmonics for metallic and metallic-colored ratios}
\label{sect:h1h2_metallic}

This section is devoted to the proof of Theorem~\ref{thm:corners}
and the justification of Numerical Result~\ref{numres:corners2},
providing a more accurate description of the functions
$h_1(\eps)$ and~$h_2(\eps)$ for the cases of
1-periodic and 2-periodic continued fractions,
i.e.~for metallic and metallic-colored ratios $\Omega$,
introduced in Section~\ref{sect:main},
using some arithmetic results from Section~\ref{sect:metallic}.
We emphasize the different behavior of the two functions in each case.

The main issue is to discuss whether the first dominant harmonic $S_1(\eps)$
and (eventually) the second one $S_2(\eps)$ are given, for any $\eps$,
by primary resonances: \ $S_i(\eps)=\bS_i(\eps)$.
\ If so, the function $h_1(\eps)$ and (eventually) the function $h_2(\eps)$
coincide with the functions $\bh_i(\eps)$ introduced in~(\ref{eq:bh1bh2}),
whose description is very simple (as in Figure~\ref{fig:omega3}).
Otherwise, the dominant harmonics are given by secondary resonances
at least for some intervals of $\eps$,
which leads to a more complicated function $h_2(\eps)$
(than $\bh_2(\eps)$, as in Figure~\ref{fig:omega13-23}\subref{fig:omega13}),
or both complicated functions $h_1(\eps)$ and $h_2(\eps)$
(as in Figures~\ref{fig:omega13-23}\subref{fig:omega23}
and~\ref{fig:omega122}).

The proof of such results requires a careful analysis of the role of
secondary resonances, is carried out rigorously
for the statements of Theorem~\ref{thm:corners}.
Instead, for the statements
of Numerical Result~\ref{numres:corners2},
we provide an evidence after having checked them numerically
for a large number of frequency ratios.

\proofjustifof{Theorem~\ref{thm:corners}(a)}
  {Numerical Result~\ref{numres:corners2}(a)}
Both of these results concern the behavior of $h_1(\eps)$
for a metallic ratio or a golden-colored ratio.
In the first case we provide a rigorous proof,
and in the second case the result relies
on Numerical Result~\ref{numres:metcol2},
which has been validated numerically for a large number of cases.

It will be enough to show the lower bound
\beq\label{eq:B0J1}
  \sqrt{B_0}>J_1\,,
\eeq
where $B_0$ is the separation between the primary and
the essential secondary resonances
(recall the definition~(\ref{eq:defB0})).
This lower bound ensures that the
most dominant harmonic is found, for all $\eps$, among the primary resonances:
$S_1(\eps)=\bS_1(\eps)$, and hence $h_1(\eps)=\bh_1(\eps)$
(such facts are reflected in Figures~\ref{fig:omega3}\subref{fig:omega3b}
and~\ref{fig:omega13-23}\subref{fig:omega13}),
which completes the proof, in view of the properties
of the function $\bh_1$ defined in~(\ref{eq:bh1bh2}).

Hence, it remains to show that the inequality~(\ref{eq:B0J1}) is fulfilled
in the two cases of a metallic and a golden-colored ratio.
Notice that, by the definition of $J_1$ in~(\ref{eq:defJ}),
we can rewrite~(\ref{eq:B0J1}) as \ $4B_0\lambda>(\lambda+1)^2$.

For a metallic ratio \ $\Omega=\ppcf{a}$, \ $a\ge1$,
\ we know from Proposition~\ref{prop:metallic} that
$B_0=5$ if $a=1$, and $B_0=a$ if $a\ge2$ (a~rigorous result).
Then, the inequality~(\ref{eq:B0J1}) is easily checked
using the quadratic equation~(\ref{eq:Tmetallic2}).

On the other hand, for a golden-colored ratio \ $\Omega=\ppcf{1,b}$, \ $b\ge2$,
\ by Numerical Result~\ref{numres:metcol2} we have $B_0=b$
(which has been established for $2\le b\le10^6$).
Then, it is easy to check the inequality~(\ref{eq:B0J1})
using in this case the quadratic equation~(\ref{eq:Tmetcol2}).
\qed

\proofof{Theorem~\ref{thm:corners}(b)}
Let us consider a golden-colored ratio \ $\Omega=\ppcf{1,b}$, \ $b\ge2$.
\ We know from Proposition~\ref{prop:metcol}
that the primary resonances $s_0(n)=s(\wh q,n)$
and the main secondary resonances $s_1(n)=s(\ol q,n)$ are generated,
respectively, by $\wh q=1$ and $\ol q=b+1$,
or, equivalently, by the vectors $v(1)=(-1,1)$ and $v(2)=(-b,b+1)$.
To study the relative position of the
graphs of the functions $g^*_{s_1(n)}$
with respect to the functions $\bg_n=g^*_{s_0(n)}$,
we compute:
\[
  \frac{\eps^*_{s_1(n)}}{\beps_n}
  =\frac{\p{\tl\gamma^*_{b+1}}^2K_1^{\,4}}{K_{b+1}^{\,4}}
  =\frac{b^2(\lambda-1)^4}{(b\lambda)^4}
  =\frac1{\lambda^2}\,,
\]
where we have used the quadratic equation~(\ref{eq:Tmetcol2})
and the fact that, by~(\ref{eq:TU4}),
\[
  \frac{K_{b+1}}{K_1}
  =\frac{z_{b+1}}{z_1}
  =\frac{b(b+1)+b\Omega}{b+\Omega}
  =\frac{b\lambda}{\lambda-1}\,.
\]
Hence, we have seen that \ $\eps^*_{s_1(n)}=\beps_n/\lambda^2=\peps_{n-1}$,
\ i.e.~the geometric means introduced in~(\ref{eq:In})
(see also Figure~\ref{fig:omega13-23}\subref{fig:omega13}).

Now, let us check that
\beq\label{eq:B0J1J2}
  \sqrt{B_0}\,J_1<J_2.
\eeq
We know from Proposition~\ref{prop:metcol}
that $B_0\le b$ (a rigorous result).
Then, it is enough to see that
\ $b\lambda(\lambda+1)^2<(\lambda^2+1)^2$,
which can be easily checked using again
the quadratic equation~(\ref{eq:Tmetcol2}).

Notice that \ $g_{s_1(n)}(\beps_n)=g_{s_1(n+1)}(\beps_n)=\dsqrt{B_0}\,J_1$.
\ We deduce from~(\ref{eq:B0J1J2}) that,
for some interval around $\beps_n$,
the second dominant harmonic $S_2$ is not the primary resonance
$\bS_2=s_0(n-1)$ (if $\eps>\beps_n$) or $\bS_2=s_0(n+1)$ (if $\eps<\beps_n$),
since at least a main secondary resonance is more dominant:
$s_1(n)$ (if $\eps>\beps_n$)
or $s_1(n+1)$ (if $\eps<\beps_n$)
(eventually, another secondary resonance
could also be the second dominant harmonic~$S_2$).
This implies that at least 3~changes in the second dominant harmonic take place
in a given period, and hence the function $h_2$ has at least 3~corners.
We also deduce that the maximum value of the function $h_2$ is $<J_2$,
since the value $J_2$ can only be reached, at the points $\beps_n$,
if a primary resonance is the second dominant there
(again, see Figure~\ref{fig:omega13-23}\subref{fig:omega13}
for an illustration, where $h_2$ has 4~corners in each period).

Concerning the minimum of the function $h_2$, it is always reached
at the points $\peps_n$ by a primary resonance:
\ $\bg_n(\peps_n)=\bg_{n-1}(\peps_n)=J_1$,
\ since by the inequality~(\ref{eq:B0J1}) all secondary resonances
take greater values at $\peps_n$.
\qed

\proofof{Theorem~\ref{thm:corners}(c)}
Now, we consider a metallic-colored but not golden-colored ratio,
\ $\Omega=\ppcf{a,b}$, \ $2\le a<b$.
\ In this case, we know from Proposition~\ref{prop:metcol}
that the primary and main secondary resonances, $s_0(n)$ and $s_1(n)$,
are generated, respectively, by $\wh q=a$ and $\ol q=1$,
or, equivalently, by the vectors $v(1)=(-1,a)$ and $v(0)=(0,1)$.
As in part~(c), we study the relative position of the functions $g^*_{s_1(n)}$
with respect to $\bg_n$, by computing:
\[
  \frac{\eps^*_{s_1(n)}}{\beps_n}
  =\frac{\p{\tl\gamma^*_1}^2K_a^{\,2}}{K_1^{\,4}}
  =\frac{(b/a)^2(\lambda-1)^4}{b^4}
  =\lambda^2,
\]
where we have used the quadratic equation~(\ref{eq:Tmetcol2})
and the fact that
\[
  \frac{K_a}{K_1}
  =\frac{z_a}{z_1}
  =\frac{ab+a\Omega}b
  =\frac{\lambda-1}b\,.
\]
Hence, we have seen that \ $\eps^*_{s_1(n)}=\lambda^2\beps_n=\peps_n$
(see also Figure~\ref{fig:omega13-23}\subref{fig:omega23}).

Next, we show that, instead of~(\ref{eq:B0J1}),
we have
\[
  \sqrt{B_0}<J_1
\]
or, equivalently,
\ $4B_0\lambda<(\lambda+1)^2$.
\ We know from Proposition~\ref{prop:metcol}
that $B_0\le b/a$ (a rigorous result).
Then, it is enough to see that
\ $4b\lambda<a(\lambda+1)^2$,
which can be checked using again the quadratic equation~(\ref{eq:Tmetcol2}).

We deduce that, for some interval around $\peps_n$,
the most dominant harmonic is not a primary resonance: $S_1\ne\bS_1$,
since at least the secondary resonance $s_1(n)$ is more dominant.
This implies that at least 2~changes
in the dominance take place in a given period, and hence the function $h_1$
has at least 2~corners.
We also deduce that the maximum value of the function $h_1$ is $<J_1$,
since the value $J_1$ can only be reached at the points $\peps_n$,
provided a primary resonance is the most dominant there
(again, see Figure~\ref{fig:omega13-23}\subref{fig:omega23}
for an illustration).
\qed

\justifof{Numerical Result~\ref{numres:corners2}(b)}
We consider a metallic ratio \ $\Omega=\ppcf{a}$, \ $a\ge1$,
\ and we are going to show that, for any~$\eps$,
the second dominant harmonic is also a primary resonance:
$S_2(\eps)=\bS_2(\eps)$, and hence $h_2(\eps)=\bh_2(\eps)$
(see Figure~\ref{fig:omega3} for an illustration).
Then, it is enough to use the simple properties
of the function $\bh_2$ defined in~(\ref{eq:bh1bh2}).

Thus, we have to check that a secondary resonance cannot be the second dominant
harmonic in any interval of~$\eps$. By the periodicity, we can restrict
ourselves to primitive vectors: $s(q,0)=k^0(q)$.
The function $g^*_{s(q,0)}$ reaches its minimum at the point $\eps^*_{s(q,0)}$,
belonging for some $n=n(q)$ to one of the intervals $\I_n=\I^+_n\cup\I^-_n$
(see the definition~(\ref{eq:In})).
Assume for instance that \ $\eps^*_{s(q,0)}\in\I^+_n=\pq{\peps_{n+1},\beps_n}$.
\ In this interval, the two dominant harmonics among the primary
resonances are $\bS_1=s_0(n)=S_1$ and $\bS_2=s_0(n+1)$.
By the inequality~(\ref{eq:B0J1}),
the second dominant harmonic among all resonances is $S_2=\bS_2$,
at least for $\eps\in\I^+_n$ close enough to $\peps_{n+1}$,
and we have to check that this is also true on the whole interval~$\I^+_n$.
Otherwise, assume that $S_2=s(q,0)$ (a secondary resonance)
for some values $\eps\in\I^+_n$ (far from~$\peps_{n+1}$).
Then, there would be an intersection
between the graphs of $\bg_{n+1}$ and $g^*_{s(q,0)}$
in the interval $\I^+_{n+1}$ and,
in view of the uniqueness given by Lemma~\ref{lm:ZW},
we would have $g^*_{s(q,0)}(\beps_n)<\bg_{n+1}(\beps_n)=J_2$.
A~symmetric discussion can be done for the case $\eps^*_{s(q,0)}\in\I^-_n$.

By the above considerations, we have to check that,
for any essential primitive $q$, and denoting $n=n(q)$ as above,
we have the lower bound
\beq\label{eq:proof2b}
  g^*_{s(q,0)}(\beps_n)\ge J_2.
\eeq
Since the minimal value of the function $g^*_{s(q,0)}$
is $\p{\tl\gamma^*_q}^{1/2}$, it is enough to consider the
essential primitives such that $\p{\tl\gamma^*_q}^{1/2}<J_2$
(by the lower bound~(\ref{eq:lowerbound}),
there is a finite number of such primitives).
We have carried out a numerical verification of~(\ref{eq:proof2b})
for all metallic ratios with $1\le a\le10^4$.

Finally, we have to justify the the statement concerning
the number of zeros $\theta^*$ of the splitting function $\M(\theta)$,
for any $\eps$ except for a small neighborhood
of the transition values $\weps$.
Notice that, since the second dominant harmonic changes from
$s_0(n-1)$ to $s_0(n+1)$ as $\eps$ goes from $\I^-_n$ to $\I^+_n$,
the transition values are $\weps=\beps_n$.
As mentioned above, the dominant harmonics
$S_1=\bS_1=s_0(n)$ and $S_2=\bS_2=s_0(n\pm1)$ are two consecutive
resonant convergents, and hence we get $\kappa=1$ in~(\ref{eq:kappa}).
As explained in Section~\ref{sect:transv}, this implies directly
that the number of zeros of $\M(\theta)$ is exactly~$4\kappa=4$.
\qed

\bremark
The exact equality in~(\ref{eq:proof2b}) holds true for some primitives $q$,
as one can see in Figure~\ref{fig:omega3} for the concrete case of
the bronze ratio.
\eremark

To end this section, we point out that the difficulty in
this kind of results comes
from the fact that we are considering infinite families of frequency ratios.
But for a particular frequency ratio, it is always possible
to obtain rigorous results, after a finite analysis
of primary and secondary resonances.
As an example, for $\Omega=\ppcf{1,2,2}$
it is not hard to see that the functions $h_1(\eps)$ and $h_2(\eps)$
have 4 and 12~corners in each period, respectively,
as well as to give explicit formulas for their minimum and maximum values
(this is illustraded in Figure~\ref{fig:omega122}).

\section{Justification of the asymptotic estimates}\label{sect:technical}

\subsection{Approximation of the splitting potential by its dominant harmonics}

The last part of this paper is devoted to the proof of Theorem~\ref{thm:main},
which gives exponentially small asymptotic estimates
of the maximal splitting distance
and the transversality of the splitting.
We start with describing our approach in a few words.

Notice that Theorem~\ref{thm:main} is stated in terms
of the Fourier coefficients of the splitting function $\M=\nabla\Lc$
introduced in~(\ref{eq:defM}).
We write, for the splitting potential and function,
\beq\label{eq:Lk}
  \Lc(\theta)=\sum_{k\in\Zc\setminus\pp0}\Lc_k\,\cos(\scprod k\theta-\tau_k),
  \quad
  \M(\theta)=-\!\!\sum_{k\in\Zc\setminus\pp0}\M_k\,\sin(\scprod k\theta-\tau_k),
\eeq
with scalar positive coefficients $\Lc_k$,
and vector coefficients
\beq\label{eq:Mk}
  \M_k=k\,\Lc_k.
\eeq
Although the Melnikov approximation~(\ref{eq:melniapproxM})
is in principle valid for real $\theta$, it is standard to see that it
can be extended to a complex strip of suitable width
(see for instance \cite{DelshamsGS04}), from which one gets
upper bounds for $\abs{\Lc_k-\mu L_k}$ and $\abs{\tau_k-\sigma_k}$,
which imply the asymptotic estimates given below in Lemma~\ref{lm:dominantsL},
ensuring that the most dominant harmonics of the Melnikov potential $L(\theta)$
are also dominant for the splitting potential~$\Lc(\theta)$.
The asymptotic estimates for the maximal splitting distance
and the transversality, given in Theorem~\ref{thm:main},
are determined from a few (one or two) dominant harmonics of the potential.
Thus, we consider approximations on $\Lc(\theta)$ given
by such dominant harmonics, together with estimates
of the sum of all other harmonics,
which show that they are dominated by the most dominant~ones.

For the proof of part~(a) of the theorem,
that provides an asymptotic estimate for the maximal splitting distance,
it will be enough to consider the approximation given by
the first dominant harmonic. Thus, we write
\beq\label{eq:L1}
  \Lc(\theta)=\Lc^{(1)}(\theta)+\F^{(2)}(\theta),
  \qquad
  \Lc^{(1)}(\theta):=\Lc_{S_1}\cos(\scprod{S_1}\theta-\tau_{S_1}),
\eeq
and we give below, in Lemma~\ref{lm:dominantsL},
an estimate of the sum of all harmonics in the remainder $\F^{(2)}(\theta)$.
This ensures that
the maximal splitting distance can be approximated by the size
of the coefficient of the dominant harmonic $S_1=S_1(\eps)$
(see the proof of Theorem~\ref{thm:main}(a) below).

On the other hand, for the proof of parts~(b) and~(c)
of Theorem~\ref{thm:main},
which concern the transversality of the splitting,
we need to detect simple zeros of $\M(\theta)$.
This is not possible with the approximation~(\ref{eq:L1})
given by only one harmonic, and we need to consider at least two harmonics.
Recalling that, by~(\ref{eq:Mk}), each (vector) harmonic $\M_k$ of $\M(\theta)$
lies in the direction of the integer vector $k$,
we consider the two most dominant \emph{essential} harmonics,
given by \emph{linearly independent} quasi-resonances
$S_1=S_1(\eps)$ and~$S_2=S_2(\eps)$
(recall the definition of essential quasi-resonances
at the beginning of Section~\ref{sect:resonantseq}).
In fact, some non-essential harmonics $c\,S_1$, $c=2,\ldots,m$,
can be (eventually) more dominant than~$S_2$.
In order to show, in Section~\ref{sect:transv}, that
such non-essential harmonics have no effect on the transversality,
we consider them separately, with specific upper bounds.
We define the \emph{index of non-essentiality} $m=m(\eps)\ge1$
as the integer satisfying
\[
  g^*_{S_1}(\eps)<\cdots<g^*_{mS_1}(\eps)\le g^*_{S_2}(\eps)<g^*_{(m+1)S_1}.
\]
Recall from~(\ref{eq:nonessential1}) that \ $g^*_{cS_1}=c\,g^*_{S_1}$.
\ It is clear that $m=1$ if and only if the two most dominant harmonics
are the non-essential ones $S_1$ and $S_2$.
For instance, we see from Figure~\ref{fig:omega13-23}\subref{fig:omega13}
that, for the case $\Omega=\ppcf{1,3}$, we have $m=2$ for $\eps$
belonging to some intervals, and $m=1$ for the remaining values of~$\eps$.

Now, we write
\bea
  \label{eq:L2}
  &&\Lc(\theta)=\Lc^{(2)}(\theta)+\F^{(\wh2)}(\theta)+\F^{(3)}(\theta),
\\
  \nonumber
  &&\Lc^{(2)}(\theta)
  :=\Lc_{S_1}\cos(\scprod{S_1}\theta-\tau_{S_1})
    +\Lc_{S_2}\cos(\scprod{S_2}\theta-\tau_{S_2}),
\\
  \nonumber
  &&\F^{(\wh2)}(\theta)
  :=\sum_{c=2}^m\Lc_{cS_1}\cos(c\scprod{S_1}\theta-\tau_{cS_1}),
\eea
and $\F^{(3)}(\theta)$ containing all harmonics not in
$\Lc^{(2)}(\theta)$ or $\F^{(\wh2)}(\theta)$
(of course, we consider $\F^{(\wh2)}=0$ if $m=1$).
Then, suitable estimates of the harmonics
in $\F^{(\wh2)}(\theta)$ and $\F^{(3)}(\theta)$,
allow us to establish the existence
of simple zeros $\theta^*$ of $\M(\theta)$,
together with an asymptotic estimate for
the minimal eigenvalue of $\Df\M(\theta^*)$,
which can be taken as a measure for the transversality of the splitting
(see Section~\ref{sect:transv}).
However, we have to exclude some intervals where
the second and third essential dominant harmonics
are of the same magnitude and the approximation~(\ref{eq:L2})
does not ensure transversality.
Such intervals are small neighborhoods of the \emph{transition values} $\weps$,
where a change in the second dominant harmonic takes place.
Such transition values can be defined as the values where
\beq\label{eq:transition}
  h_2(\weps)=h_3(\weps).
\eeq

We will use the following lemma, analogous
to the one established in \cite{DelshamsG03,DelshamsG04},
providing an asymptotic estimate
for the dominant harmonics $\Lc_{S_1}$ and $\Lc_{S_2}$
(and an upper bound for the difference of their phases $\tau_{S_i}$
with respect to the original ones~$\sigma_{S_i}$),
as well as an estimate for the sum of all the harmonics
in the remainders $\F^{(i)}$, $i=2,\wh2,3$,
appearing in~(\ref{eq:L1}) and~(\ref{eq:L2}).
To unify the notation, we write
\ $\ds\F^{(i)}=\sum_{k\in\Zc_i}(\cdots)$,
\ defining the sets of indices
\[
  \Zc_2=\Zc\setminus\{0,S_1\},
  \qquad
  \Zc_{\wh2}=\{2S_1,\cdots,m S_1\},
  \qquad
  \Zc_3=\Zc\setminus(\{0,S_1,S_2\}\cup\Zc_{\wh2}).
\]
The estimate for each sum is given,
due to the exponential smallness of the harmonics,
in terms of the dominant harmonic in each set $\Zc_i$,
that we denote as $\wtS_i=\wtS_i(\eps)$, $i=2,\wh2,3$.
Notice that, for $i=\wh2$, the dominant harmonic
is clearly $\wtS_{\wh2}=2S_1$ (non-essential)
and, for $i=2,3$, the dominant harmonic can be either
$\wtS_i=S_i$ (essential) or $\wtS_i=2S_{i-1}$ (non-essential).
With this in mind, we introduce the functions
\beq\label{eq:wth}
  \begin{array}{ll}
    \wth_{\wh2}(\eps):=2h_1(\eps)=g^*_{\wtS_{\wh2}}(\eps).
  \\[10pt]
    \wth_i(\eps):=\min(h_i(\eps),\,2h_{i-1}(\eps))=g^*_{\wtS_i}(\eps),
    \quad i=2,3,
  \end{array}
\eeq
We stress that, in the three cases, the function $\wth_i(\eps)$ is given
by the minimum of the values $g^*_k(\eps)$, with $k$ belonging to the
corresponding set of indices:
\ $\wth_i(\eps)=\ds\min_{k\in\Zc_i}g^*_k(\eps)$, \ $i=2,\wh2,3$.
\ Comparing with the functions $h_i(\eps)$ defined in~(\ref{eq:h1h2}),
we see that non-essential harmonics are also taken into account
in the definition of $\wth_i(\eps)$.
Notice also that the equality~(\ref{eq:transition}) characterizing
the transition values can be rewritten as $h_2(\weps)=\wth_3(\weps)$.

Recall that the coefficients $\Lc_k$, introduced in~(\ref{eq:Lk}),
are all positive.
In fact, we are not directly interested in the
splitting potential $\Lc(\theta)$, but rather in some of its derivatives
(such as $\M(\theta)$, $\Df\M(\theta)$).
The constant $C_0$ in the exponentials is the one
defined in~(\ref{eq:beta_gk}).
On the other hand, recall that the meaning of the notations `$\sim$'
and `$\preceq$' has been introduced at the end of Section~\ref{sect:main}.

\begin{lemma}\label{lm:dominantsL}
For $\eps$ small enough and $\mu=\eps^r$ with $r>3$, one has:
\btm
\item[\rm(a)]
$\ds\Lc_{S_i}
 \sim\mu\,L_{S_i}
 \sim\frac{\mu}{\eps^{1/4}}
   \,\exp\pp{-\frac{C_0h_i(\eps)}{\eps^{1/4}}}$,
\quad
$\abs{\tau_{S_i}-\sigma_{S_i}}\preceq\dfrac\mu{\eps^3}$\,,
\quad $i=1,2$;
\vspace{4pt}
\item[\rm(b)]
$\ds\sum_{k\in\Zc_i}\Lc_k
\sim\frac1{\eps^{1/4}}\,\Lc_{\wtS_i}
\sim\frac\mu{\eps^{1/4}}
  \,\exp\pp{-\frac{C_0\wth_i(\eps)}{\eps^{1/4}}}$,
\quad $i=2,\wh2,3$.
\etm
\end{lemma}

\sketchproof
We only give the main ideas of the proof, since it is
similar to analogous results in \cite[Lemmas~4 and~5]{DelshamsG04}
and \cite[Lemma~3]{DelshamsG03}.
At first order in $\mu$,
the coefficients of the splitting potential can be approximated,
neglecting the error term in the Melnikov
approximation~(\ref{eq:melniapproxM}),
by the coefficients of the Melnikov potential,
given in~(\ref{eq:alphabeta}):
\ $\Lc_k\sim\mu L_k=\mu\alpha_k\,\ee^{- \beta_k}$.
\ As mentioned in Section~\ref{sect:gn},
the main behavior of the coefficients~$L_k(\eps)$ is given by
the exponents $\beta_k(\eps)$, which have been written in~(\ref{eq:beta_gk})
in terms of the functions $g_k(\eps)$.
In particular, the coefficients $L_{S_i}$
associated to the two essential dominant
harmonics $k=S_i(\eps)$, $i=1,2$,
can be expressed in terms of the functions $h_i(\eps)$
introduced in~(\ref{eq:h1h2}).
In this way, we obtain an estimate for the factor $\ee^{-\beta_{S_i}}$,
which provides the exponential factor in~(a).

We also consider the factor $\alpha_k$,
with $k=S_i(\eps)$. Recalling from~(\ref{eq:estimS}) that
\ $\abs{S_i}\sim\eps^{-1/4}$,
we get from~(\ref{eq:alphabeta}) that
\ $\alpha_{S_i}\sim\eps^{-1/4}$,
\ which provides the polynomial factor in part~(a).

The estimate obtained is valid for the dominant coefficient
of the Melnikov potential $L(\theta)$.
To complete the proof of part~(a),
one has to show that an analogous estimate is also valid
for the splitting potential $\Lc(\theta)$,
i.e.~when the error term in the Poincar\'e--Melnikov
approximation~(\ref{eq:melniapproxM}) is not neglected.
This requires to obtain an upper bounds
(provided in \cite[Th.~10]{DelshamsGS04})
for the corresponding coefficient of the error term in~(\ref{eq:melniapproxM})
and show that, in our singular case $\mu=\eps^r$, it is also exponentially
small and dominated by the main term in the approximation.
This can be worked out straightforwardly as in \cite[Lemma~5]{DelshamsG04}
(where the case of the golden number was considered),
so we omit the details here.

The proof of part~(b) is carried out in similar terms.
For the dominant harmonic $k=\wtS_i$ inside each set
$\Zc_i$, $i=2,\wh2,3$,
we get $\abs{\wtS_i}\sim\eps^{-1/4}$ as in~(\ref{eq:estimS}),
and an exponentially small estimate
for $\Lc_{\wtS_i}$ with the function $\wth_i(\eps)$
defined in~(\ref{eq:wth}).
Such estimates are also valid if one considers the whole sum in~(b),
since for any given $\eps$ the terms of this sum can be bounded
by a geometric series and, hence, it can be estimated by its dominant term
(see \cite[Lemma~4]{DelshamsG04} for more details).
\qed

In regard to the proof of Theorem~\ref{thm:main}(a,c),
we need to measure the size of each perturbation $\F^{(i)}(\theta)$
in~(\ref{eq:L1}--\ref{eq:L2})
with respect to the coefficients of the approximations $\K^{(j)}(\theta)$.
Since by Lemma~\ref{lm:dominantsL} the size of $\F^{(i)}(\theta)$
is given by the size of its dominant harmonic,
we introduce the following small parameters:
\beq\label{eq:defetaij}
  \eta_{i,j}:=\frac{\Lc_{\wtS_i}}{\Lc_{S_j}}
  \sim\exp\pp{-\frac{C_0(\wth_i(\eps)-h_j(\eps))}{\eps^{1/4}}},
  \quad
  (i,j)=(2,1),(\wh2,1),(3,1),(3,2),
\eeq
as a measure of the perturbations $\F^{(i)}$ in~(\ref{eq:L1}--\ref{eq:L2}),
relatively to the size of the essential dominant coefficients $\Lc_{S_j}$
(we~consider $\eta_{\wh2,1}=0$ if the index of non-essentiality is $m=1$).
Although we define the parameters $\eta_{i,j}$
in terms of the coefficients of~$\Lc(\theta)$,
we can also define them from the coefficients of its derivatives,
such as the splitting function~$\M(\theta)$,
in view of~(\ref{eq:Mk}) and the fact that the respective factors
have the same magnitude: $\abs{S_j}\sim\abs{\wtS_i}$.

The parameters $\eta_{i,j}$ always exponentially small in $\eps$,
provided we exclude some small neighborhoods
where $\Lc_{S_j}$ and $\Lc_{\wtS_i}$ can have the same magnitude.
For instance, we have $\eta_{3,2}$ exponentially small if $\eps$
is not very close to the transition values~(\ref{eq:transition}),
at which the second and third dominant harmonics have the same magnitude.
Analogously, the parameter $\eta_{2,1}$ is exponentially small
excluding neighborhoods where the first and second dominant harmonics
have the same magnitude.

\proofof{Theorem~\ref{thm:main}(a)}
Applying Lemma~\ref{lm:dominantsL}, we see that
the coefficient of the dominant harmonic of the splitting function $\M(\theta)$
is greater than the sum of all other harmonics.
More precisely, we have for $\eps\to0$ the estimate
\beq\label{eq:Mdom}
  \max_{\theta\in\T^2}\abs{\M(\theta)}
  =\abs{\M_{S_1}}(1+\Ord(\eta_{2,1}))
  \sim \abs{\M_{S_1}}
  \sim\abs{S_1}\Lc_{S_1},
\eeq
which implies the result,
using the asymptotic estimate~(\ref{eq:estimS}) for $\abs{S_1}$,
and the asymptotic estimate for $\abs{\M_{S_1}}$, in terms of $h_1(\eps)$,
deduced from Lemma~\ref{lm:dominantsL}(a).

We point out that the previous argument does not apply directly
when $\eps$ is close to a value where $h_1$ and $h_2$ coincide,
i.e.~the first and second dominant harmonics have the same magnitude
(for instance, for a metallic ratio $\Omega$ this occurs near the values
$\peps_n$, see~(\ref{eq:In}) and Figure~\ref{fig:omega3}\subref{fig:omega3b}).
Eventually, more than two harmonics
(but a finite number, according to the arguments given in Lemma~\ref{lm:ZW})
might also have the same magnitude and become dominant.
In such cases, the parameter $\eta_{2,1}$ is not exponentially small,
but we can replace the main term in~(\ref{eq:Mdom}) by a finite number
of terms, plus an exponentially small perturbation, and by the properties of
Fourier expansions the maximum value of~$\abs{\M(\theta)}$ can be
compared to any of its dominant harmonics.
\qed

\subsection{Nondegenerate critical points and transversality}
\label{sect:transv}

This section is devoted to the study of the transversality of the
homoclinic orbits for values of the perturbation parameter $\eps$,
not very close to the transition values $\weps$
defined in~(\ref{eq:transition}).
For such values of $\eps$ we show that, under suitable conditions,
the splitting potential $\Lc(\theta)$ has $4\kappa$ nondegenerate
critical points for some integer $\kappa\ge1$,
i.e.~the splitting function $\M(\theta)$ has $4\kappa$ simple zeros,
which give rise to $4\kappa$ transverse homoclinic orbits.
Such critical points are easily detected
in the approximation $\Lc^{(2)}(\theta)$
introduced in~(\ref{eq:L2}), given by the 2~essential dominant harmonics,
and using the estimates for $\F^{(\wh2)}(\theta)$ and $\F^{(3)}(\theta)$
given in Lemma~\ref{lm:dominantsL} we can prove the persistence
of the critical points in the whole function $\Lc(\theta)$.

In fact, we make the computations easier by
performing a linear change on $\T^2$,
taking $\Lc^{(2)}(\theta)$ to a very simple form.
As in \cite{DelshamsG03,DelshamsG04},
we introduce the variables
\beq\label{eq:psichange}
  \psi_1=\scprod{S_1}\theta-\tau_{S_1},
  \qquad
  \psi_2=\scprod{S_2}\theta-\tau_{S_2}.
\eeq
This change of variables is valid for $\eps$ in the interval between two
consecutive transition values, in which we have
two concrete essential dominant harmonics $S_1(\eps)$ and $S_2(\eps)$,
which remain constant in this interval.
In the new variables, the functions
$\Lc$, $\Lc^{(2)}$, $\F^{(\wh2)}$, $\F^{(3)}$ in~(\ref{eq:L2}) become
\beq
  \label{eq:K2}
  \K(\psi)=\K^{(2)}(\psi)+\G^{(\wh2)}(\psi_1)+\G^{(3)}(\psi)
\eeq
where, in particular, we have
\beq\label{eq:K2b}
  \K^{(2)}(\psi)=\Lc_{S_1}\cos\psi_1+\Lc_{S_2}\cos\psi_2.
\eeq

It is clear that $\K^{(2)}$ has the 4~critical points, all nondegenerate:
\ $\psi^{*,0}:=(0,0),(\pi,0),(0,\pi),(\pi,\pi)$
\ (one maximum, two saddles and one minimum, respectively).
Regarding $\K(\psi)$ as a perturbation of $\K^{(2)}(\psi)$,
we are going to show that
it also has 4~critical points $\psi^*$, all nondegenerate,
which are close to the critical points $\psi^{*,0}$ of $\K^{(2)}(\psi)$.
We point out that, in general, the change~(\ref{eq:psichange})
is not one-to-one on $\T^2$, but rather ``$\kappa$-to-one'', where
\beq\label{eq:kappa}
  \kappa=\kappa(\eps):=\abs{\det(S_1,S_2)}.
\eeq
Hence, the number of critical points of $\Lc(\theta)$ is $4\kappa$.
It is not hard to show that $\kappa(\eps)$
is $4\ln\lambda$-periodic in $\ln\eps$.
Moreover, it is ``piecewise-constant'' with (eventual)
jump discontinuities when changes in the dominant harmonics take place.

\bremark
For a metallic ratio \ $\Omega=\ppcf{a}$, \ we know from
Numerical Result~\ref{thm:corners}(b) that $\kappa=1$
(a result checked numerically for $1\le a\le10^4$),
and hence there are exactly 4~transverse homoclinic orbits, for any $\eps$
small enough (excluding a neighborhood of the transition values $\weps$).
Although for other frequency ratios it is possible, in principle,
to have $\kappa\ge2$, we have obtained $\kappa=1$
for all the cases we have explored.
\eremark

To establish the persistence of the critical points,
we are going to use the following lemma, whose proof is a simple application
of the 2-dimensional fixed point theorem and is omitted here.

\begin{lemma}\label{lm:FPsin2D}
If \ $f_1,f_2:\T^2\longrightarrow\R$ \ are differentiable and satisfy
\[
  f_i^{\,2}+\p{\abs{\pderiv{f_i}{\psi_1}}+\abs{\pderiv{f_i}{\psi_2}}}^2<1,
  \qquad
  i=1,2,
\]
then the system of equations
\beq\label{eq:FPsin2D}
  \sin\psi_1=f_1(\psi),
  \qquad
  \sin\psi_2=f_2(\psi)
\eeq
has exactly 4 solutions $\psi^*$, which are simple.
Furthermore, if \ $f_1(\psi),f_2(\psi)=\Ord(\eta)$ \ for any $\psi\in\T^2$,
with $\eta$ sufficiently small, then the solutions of the system satisfy
\ $\psi^*=\psi^{*,0}+\Ord(\eta)$,
\ with \ $\psi^{*,0}=(0,0),(0,\pi),(\pi,0),(\pi,\pi)$.
\end{lemma}

In order to apply this lemma, we introduce the following
perturbation parameter, using the parameters $\eta_{i,j}$
defined in~(\ref{eq:defetaij}),
\beq\label{eq:defeta}
  \eta:=\max(\eta_{\wh2,1},\eta_{3,1},\eta_{3,2}).
\eeq

\begin{lemma}\label{lm:critpK}
The function $\K(\psi)$ has exactly 4~critical points, all nondegenerate:
\beq\label{eq:psipoints}
  \psi^*=\psi^{*,0}+\Ord(\eta),
  \qquad
  \mbox{with }\ \psi^{*,0}=(0,0),(0,\pi),(\pi,0),(\pi,\pi).
\eeq
At each critical point, we have
\[
  \det\Df\K(\psi^*)
  =\delta^*_1\delta^*_2\,\Lc_{S_1}\,\Lc_{S_2}
    \p{1+\Ord(\eta)},
\]
where \ $\delta^*_i=\cos\psi^{*,0}_i=\pm1$, \ $i=1,2$.
\end{lemma}

\proof
We see from~(\ref{eq:K2}--\ref{eq:K2b})
that the system of equations $\nabla\K(\psi)=0$
can be written as in~(\ref{eq:FPsin2D}), with the functions
\[
  f_1(\psi)
  =\frac1{\Lc_{S_1}}\p{\deriv{\G^{(\wh2)}}{\psi_1}+\pderiv{\G^{(3)}}{\psi_1}},
  \qquad
  f_2(\psi)=\frac1{\Lc_{S_2}}\cdot\pderiv{\G^{(3)}}{\psi_2}
\]
(notice that $\G^{(\wh2)}$ does not depend on $\psi_2$).
By Lemma~\ref{lm:dominantsL}, we have $f_1(\psi),f_2(\psi)=\Ord(\eta)$,
with $\eta$ as given in~(\ref{eq:defeta}).
Hence, applying Lemma~\ref{lm:FPsin2D} we deduce the result for the
critical points of $\K(\psi)$.

We also provide, for each critical point, an asymptotic estimate
for the determinant of $\Df^2\K(\psi^*)$.
It is clear, for the perturbed critical points, that
the signs $\delta^*_i=\pm1$ become perturbed as follows:
\ $\cos\psi^*_i=\delta^*_i+\Ord(\eta^2)$.
\ Writing \ $\Df^2\K(\psi)=\symmatrix{k_{11}}{k_{12}}{k_{22}}$,
\ we have from~(\ref{eq:K2}) the approximations
\bea
  \label{eq:D2Ka}
  &&k_{11}=\pdderiv\K{\psi_1}
  =-\Lc_{S_1}\p{\cos\psi_1+\Ord(\eta_{\wh2,1},\eta_{3,1})},
\\
  \label{eq:D2Kb}
  &&k_{12}=\pdderivm\K{\psi_1}{\psi_2}
  =\Lc_{S_1}\cdot\Ord(\eta_{3,1})
  =\Lc_{S_2}\cdot\Ord(\eta_{3,2}),
\\
  \label{eq:D2Kc}
  &&k_{22}=\pdderiv\K{\psi_2}
  =-\Lc_{S_2}\p{\cos\psi_2+\Ord(\eta_{3,2})},
\eea
and we deduce the expression for the determinant of $\Df^2\K(\psi^*)$.
\qed

Now we complete the proof of part~(b) of our main theorem
by applying the inverse of the linear change~(\ref{eq:psichange})
to the critical points $\psi^*$ of $\K(\psi)$,
in order to get the critical points $\theta^*$ of
the splitting potential $\Lc(\theta)$,
i.e.~the zeros of the Melnikov function $\M(\theta)$.

\proofof{Theorem~\ref{thm:main}(b)}
Since the linear change~(\ref{eq:psichange}) is ``$\kappa$-to-one'',
with $\kappa$ as in~(\ref{eq:kappa}),
the 4~critical points $\psi^*$ of~$\K(\psi)$ give rise
to 4$\kappa$~critical points $\theta^*$ of $\Lc(\theta)$.
It is clear that such critical points are also nondegenerate,
and hence they are simple zeros of the splitting function~$\M(\theta)$.
\qed

\bremarks
\item
Recall that the vectors $S_i=S_i(\eps)$ remain constant between consecutive
transition values $\weps$ (see the proof of Proposition~\ref{prop:h1h2}).
On the other hand, by~(\ref{eq:psipoints}) the points $\psi^*$
are $\Ord(\eta)$-close to the points $\psi^{*,0}$,
where $\eta$ is exponentially small.
Hence, the points $\theta^*=\theta^*(\eps)$ remain ``nearly constant'' along
each interval of $\eps$ between consecutive transition values $\weps$,
and can ``change'' when $\eps$ goes across a value $\weps$.
\item
As a particular interesting case, we may consider the phases $\sigma_k=0$
in the perturbation~(\ref{eq:hf}). In this case, our Hamiltonian system given
by~(\ref{eq:HamiltH}--\ref{eq:hf}) is \emph{reversible}
with respect to the involution
\beq\label{eq:reversible}
  \Rc:(x,y,\varphi,I)\mapsto(-x,y,-\varphi,I)
\eeq
(indeed, its associated Hamiltonian field satisfies
the identity $X_H\circ\Rc=-\Rc\,X_H$).
We point out that reversible perturbations
have also been considered in some related papers
\cite{Gallavotti94,GallavottiGM99a,RudnevW98}.
Under the reversibility~(\ref{eq:reversible}),
the whiskers are related by the involution: $\W^\st=\Rc\,\W^\ut$.
Hence, their parameterizations in~(\ref{eq:defM})
can be chosen in such a way that
$\J^\st(\theta)=\J^\ut(-\theta)$, provided
the transverse section $x=\pi$ is considered in their definition.
This implies that the splitting function is an odd function:
\ $\M(-\theta)=-\M(\theta)$
\ (and the splitting potential $\Lc(\theta)$ is even)
and, using its periodicity, one sees that $\M(\theta)$ has,
at least, the following 4~zeros:
\ $\theta^*\ =\ (0,0),\;(\pi,0),\;(0,\pi),\;(\pi,\pi)$
\ (notice that they do not depend on $\eps$).
Although such zeros could be non-simple in principle,
the result of Theorem~\ref{thm:main}(b)
says that they are simple for any $\eps$ except for a small neighborhood
of the transition values $\weps$,
at which some bifurcations of the zeros could take place.
\eremarks

It remains to provide, for each zero $\theta^*$ of
the splitting function $\M(\theta)$,
an asymptotic estimate for the minimal eigenvalue of
the matrix $\Df^2\Lc(\theta^*)=\Df\M(\theta^*)$,
as a measure for the transversality of the splitting.

\proofof{Theorem~\ref{thm:main}(c)}
Denoting \ $D=\det\Df^2\Lc(\theta^*)$ \ and \ $T = \tr \Df^2\Lc(\theta^*)$,
\ we can present the (modulus of the) minimal eigenvalue of
$\Df^2\Lc(\theta^*)$ in the form
\[
  \abs{m^*}=\frac{2\abs D}{\abs T+\dsqrt{T^2-4D}}\sim\frac{\abs D}{\abs T}\,,
\]
where we have taken into account that \ $0\le\dsqrt{T^2-4D}\le\abs T$.
\ Thus, we need to find estimates for $\abs D$ and $\abs T$,
at the critical points $\theta^*$ of $\Lc(\theta)$ (or zeros of $\M(\theta)$).

By the linear change~(\ref{eq:psichange}),
we have \ $\Df^2\Lc(\theta^*)=A\tp\,\Df^2\K(\psi^*)\,A$,
\ where $A$ is the matrix having the vectors $S_1$ and~$S_2$ as rows.
In~(\ref{eq:kappa}), we have defined \ $\kappa=\abs{\det A}$.
Since $\kappa=\kappa(\eps)$ is picewise-constant and periodic in $\ln\eps$,
it is bounded from below and from above: $\kappa\sim1$.
Applying Lemma~\ref{lm:critpK}, we get the asymptotic estimate
\[
  \abs D
  =\Lc_{S_1}\,\Lc_{S_2}\p{1+\Ord(\eta)}
  \sim\Lc_{S_1}\,\Lc_{S_2}.
\]

On the other hand, in order to estimate $T$ we write $\Df^2\K$
as in~(\ref{eq:D2Ka}--\ref{eq:D2Kc}), and obtain
\[
  \Df^2\Lc(\theta)
  =k_{11}\,S_1\cdot S_1\tp
   +k_{12}(S_1\cdot S_2\tp+S_2\cdot S_1\tp)
   +k_{22}\,S_2\cdot S_2\tp,
\]
which implies that
\ $T=k_{11}\abs{S_1}_2^{\,2}+2k_{12}\scprod{S_1}{S_2}+k_{22}\abs{S_2}_2^{\,2}$,
\ where $\abs\cdot_2$ denotes the usual Euclidean norm
(which is equivalent to the norm $\abs\cdot=\abs\cdot_1$
mainly used in this paper).
Now we use, at the critical points $\psi^*$,
the estimates for the matrix $\Df^2\K(\psi^*)$
given in~(\ref{eq:D2Ka}--\ref{eq:D2Kc}).
We obtain \ $\abs{k_{11}}\sim\Lc_{S_1}$ \ as the main entry, and
\ $\abs{k_{12}}\sim\Lc_{S_1}\cdot\Ord(\eta_{3,1})$,
\ $\abs{k_{22}}\sim\Lc_{S_2}$.
Applying also the estimate~(\ref{eq:estimS}), we obtain
\[
  \abs T\sim\frac1{\sqrt\eps}\,\Lc_{S_1},
  \qquad
  \mbox{and hence}
  \quad
  m^*\sim\frac{\abs D}{\abs T}\sim\sqrt{\eps}\,\Lc_{S_2}.
\]
Applying the estimate for $\Lc_{S_2}$ given in Lemma~\ref{lm:dominantsL},
we obtain the desired estimate for the minimal eigenvalue.
\qed

\small

\begin{thebibliography}{GGM99b}

\bibitem[Arn64]{Arnold64}
V.I. Arnold.
\newblock Instability of dynamical systems with several degrees of freedom.
\newblock {\em Soviet Math. Dokl.}, 5(3):581--585, 1964.
\newblock (\emph{Dokl. Akad. Nauk SSSR}, 156:9--12, 1964).

\bibitem[BFGS12]{BaldomaFGS12}
I.~Baldom\'a, E.~Fontich, M.~Guardia, and T.M. Seara.
\newblock Exponentially small splitting of separatrices beyond {M}elnikov
  analysis: {R}igorous results.
\newblock {\em J. Differential Equations}, 253(12):3304--3439, 2012.

\bibitem[DG00]{DelshamsG00}
A.~Delshams and P.~Guti\'errez.
\newblock Splitting potential and the {P}oincar\'e--{M}elnikov method for
  whiskered tori in {H}amiltonian systems.
\newblock {\em J.~Nonlinear Sci.}, 10(4):433--476, \noopsort{00a}2000.

\bibitem[DG01]{DelshamsG01}
A.~Delshams and P.~Guti\'errez.
\newblock Homoclinic orbits to invariant tori in {H}amiltonian systems.
\newblock In C.K.R.T. Jones and A.I. Khibnik, editors, {\em Multiple-Time-Scale
  Dynamical Systems (Minneapolis, MN, 1997)}, volume 122 of {\em IMA Vol. Math.
  Appl.}, pages 1--27. Springer-Verlag, New York, 2001.

\bibitem[DG03]{DelshamsG03}
A.~Delshams and P.~Guti\'errez.
\newblock Exponentially small splitting of separatrices for whiskered tori in
  {H}amiltonian systems.
\newblock {\em Zap. Nauchn. Sem. S.-Peterburg. Otdel. Mat. Inst. Steklov.
  (POMI)}, 300:87--121, 2003.
\newblock (\emph{J. Math. Sci. (N.Y.)}, 128(2):2726--2746, 2005).

\bibitem[DG04]{DelshamsG04}
A.~Delshams and P.~Guti\'errez.
\newblock Exponentially small splitting for whiskered tori in {H}amiltonian
  systems: continuation of transverse homoclinic orbits.
\newblock {\em Discrete Contin. Dyn. Syst.}, 11(4):757--783, 2004.

\bibitem[DGG14a]{DelshamsGG14a}
A.~Delshams, M.~Gonchenko, and P.~Guti\'errez.
\newblock Exponentially small asymptotic estimates for the splitting of
  separatrices to whiskered tori with quadratic and cubic frequencies.
\newblock {\em Electron. Res. Announc. Math. Sci.}, 21:41--61,
  \noopsort{14a}2014.

\bibitem[DGG14b]{DelshamsGG14b}
A.~Delshams, M.~Gonchenko, and P.~Guti\'errez.
\newblock Exponentially small lower bounds for the splitting of separatrices to
  whiskered tori with frequencies of constant type.
\newblock {\em Internat. J. Bifur. Chaos Appl. Sci. Engrg.}, 24(8):1440011,
  12~pp., \noopsort{14b}2014.

\bibitem[DGG14c]{DelshamsGG14c}
A.~Delshams, M.~Gonchenko, and P.~Guti\'errez.
\newblock Continuation of the exponentially small transversality for the
  splitting of separatrices to a whiskered torus with silver ratio.
\newblock {\em Regul. Chaotic Dyn.}, 19(6):663--680, \noopsort{14c}2014.

\bibitem[DGJS97]{DelshamsGJS97}
A.~Delshams, V.G. Gelfreich, {\`A}.~Jorba, and T.M. Seara.
\newblock Exponentially small splitting of separatrices under fast
  quasiperiodic forcing.
\newblock {\em Comm. Math. Phys.}, 189:35--71, 1997.

\bibitem[DGS04]{DelshamsGS04}
A.~Delshams, P.~Guti\'errez, and T.M. Seara.
\newblock Exponentially small splitting for whiskered tori in {H}amiltonian
  systems: flow-box coordinates and upper bounds.
\newblock {\em Discrete Contin. Dyn. Syst.}, 11(4):785--826, 2004.

\bibitem[DLS06]{DelshamsLS06}
A.~Delshams, R.{ }de{ }la Llave, and T.M. Seara.
\newblock A geometric mechanism for diffusion in {H}amiltonian systems
  overcoming the large gap problem: heuristics and rigorous verification on a
  model.
\newblock {\em Mem. Amer. Math. Soc.}, 179(844), 2006.

\bibitem[DR98]{DelshamsR98}
A.~Delshams and R.~Ram\'{\i}rez{-}Ros.
\newblock Exponentially small splitting of separatrices for perturbed
  integrable standard-like maps.
\newblock {\em J.~Nonlinear Sci.}, 8(3):317--352, 1998.

\bibitem[DS92]{DelshamsS92}
A.~Delshams and T.M. Seara.
\newblock An asymptotic expression for the splitting of separatrices of the
  rapidly forced pendulum.
\newblock {\em Comm. Math. Phys.}, 150:433--463, 1992.

\bibitem[DS97]{DelshamsS97}
A.~Delshams and T.M. Seara.
\newblock Splitting of separatrices in {H}amiltonian systems with one and a
  half degrees of freedom.
\newblock {\em Math. Phys. Electron.~J.}, 3:\ paper~4, 40~pp., 1997.

\bibitem[Eli94]{Eliasson94}
L.H. Eliasson.
\newblock Biasymptotic solutions of perturbed integrable {H}amiltonian systems.
\newblock {\em Bol. Soc. Brasil. Mat. (N.S.)}, 25(1):57--76, 1994.

\bibitem[FP07]{FalconP07}
S.~Falc\'on and \'A. Plaza.
\newblock The $k$-{F}ibonacci sequence and the {P}ascal 2-triangle.
\newblock {\em Chaos Solitons Fractals}, 33(1):38--49, 2007.

\bibitem[Gal94]{Gallavotti94}
G.~Gallavotti.
\newblock Twistless {KAM} tori, quasi flat homoclinic intersections, and other
  cancellations in the perturbation series of certain completely integrable
  {H}amiltonian systems. {A}~review.
\newblock {\em Rev. Math. Phys.}, 6(3):343--411, 1994.

\bibitem[Gel97]{Gelfreich97}
V.G. Gelfreich.
\newblock Melnikov method and exponentially small splitting of separatrices.
\newblock {\em Phys. D}, 101(3-4):227--248, 1997.

\bibitem[GGM99a]{GallavottiGM99b}
G.~Gallavotti, G.~Gentile, and V.~Mastropietro.
\newblock Melnikov approximation dominance. {S}ome examples.
\newblock {\em Rev. Math. Phys.}, 11(4):451--461, 1999.

\bibitem[GGM99b]{GallavottiGM99a}
G.~Gallavotti, G.~Gentile, and V.~Mastropietro.
\newblock Separatrix splitting for systems with three time scales.
\newblock {\em Comm. Math. Phys.}, 202(1):197--236, 1999.

\bibitem[GL06]{GideaL06}
M.~Gidea and R.{ }de{ }la Llave.
\newblock Topological methods in the instability problem of {H}amiltonian
  systems.
\newblock {\em Discrete Contin. Dyn. Syst.}, 14(2):295--328, 2006.

\bibitem[GR03]{GideaR03}
M.~Gidea and C.~Robinson.
\newblock Topologically crossing heteroclinic connections to invariant tori.
\newblock {\em J.~Differential Equations}, 193(1):49--74, 2003.

\bibitem[GS12]{GuardiaS12}
M.~Guardia and T.M. Seara.
\newblock Exponentially and non-exponentially small splitting of separatrices
  for the pendulum with a fast meromorphic perturbation.
\newblock {\em Nonlinearity}, 25(5):1367--1412, 2012.

\bibitem[HMS88]{HolmesMS88}
P.~Holmes, J.E. Marsden, and J.~Scheurle.
\newblock Exponentially small splittings of separatrices with applications to
  {KAM} theory and degenerate bifurcations.
\newblock In {\em Hamiltonian Dynamical Systems (Boulder, CO, 1987)}, volume~81
  of {\em Contemp. Math.}, pages 213--244. Amer. Math. Soc., Providence, RI,
  1988.

\bibitem[Koc99]{Koch99}
H.~Koch.
\newblock A renormalization group for {H}amiltonians, with applications to
  {KAM} theory.
\newblock {\em Ergodic Theory Dynam. Systems}, 19(2):475--521, 1999.

\bibitem[Lan95]{Lang95}
S.~Lang.
\newblock {\em Introduction to Diophantine approximations}.
\newblock Springer-Verlag, New York, 2nd edition, 1995.

\bibitem[Laz03]{Lazutkin03}
V.F. Lazutkin.
\newblock Splitting of separatrices for the {C}hirikov standard map.
\newblock {\em Zap. Nauchn. Sem. S.-Peterburg. Otdel. Mat. Inst. Steklov.
  (POMI)}, 300:25--55, 2003.
\newblock (\emph{J. Math. Sci. (N.Y.)}, 128(2):2687--2705, 2005). The original
  Russian preprint appeared in 1984.

\bibitem[LMS03]{LochakMS03}
P.~Lochak, J.-P. Marco, and D.~Sauzin.
\newblock On the splitting of invariant manifolds in multidimensional
  near-integrable {H}amiltonian systems.
\newblock {\em Mem. Amer. Math. Soc.}, 163(775), 2003.

\bibitem[Loc90]{Lochak90}
P.~Lochak.
\newblock Effective speed of {A}rnold's diffusion and small denominators.
\newblock {\em Phys. Lett.~A}, 143(1-2):39--42, 1990.

\bibitem[Lop02]{Lopesd02a}
J.~Lopes{ }Dias.
\newblock Renormalization of flows on the multidimensional torus close to a
  {KT} frequency vector.
\newblock {\em Nonlinearity}, 15(3):647--664, 2002.

\bibitem[Mel63]{Melnikov63}
V.K. Melnikov.
\newblock On the stability of the center for time periodic perturbations.
\newblock {\em Trans. Moscow Math. Soc.}, 12:1--57, 1963.
\newblock (\emph{Trudy Moskov. Mat. Ob\v{s}\v{c}.}, 12:3--52, 1963).

\bibitem[Nie00]{Niederman00}
L.~Niederman.
\newblock Dynamics around simple resonant tori in nearly integrable
  {H}amiltonian systems.
\newblock {\em J.~Differential Equations},
  161(1):\discretionary{}{}{}\mbox{1--41}, 2000.

\bibitem[Poi90]{Poincare90}
H.~Poincar\'e.
\newblock Sur le probl\`eme des trois corps et les \'equations de la dynamique.
\newblock {\em Acta Math.}, 13:1--270, 1890.

\bibitem[PT00]{ProninT00}
A.~Pronin and D.V. Treschev.
\newblock Continuous averaging in multi-frequency slow--fast systems.
\newblock {\em Regul. Chaotic Dyn.}, 5(2):157--170, 2000.

\bibitem[RW98]{RudnevW98}
M.~Rudnev and S.~Wiggins.
\newblock Existence of exponentially small separatrix splittings and homoclinic
  connections between whiskered tori in weakly hyperbolic near-integrable
  {H}amiltonian systems.
\newblock {\em Phys.~D}, 114(1-2):3--80, 1998.

\bibitem[RW00]{RudnevW00}
M.~Rudnev and S.~Wiggins.
\newblock On a homoclinic splitting problem.
\newblock {\em Regul. Chaotic Dyn.}, 5(2):227--242, 2000.

\bibitem[Sau01]{Sauzin01}
D.~Sauzin.
\newblock A new method for measuring the splitting of invariant manifolds.
\newblock {\em Ann. Sci. \'Ecole Norm. Sup.~(4)}, 34(2):159--221, 2001.

\bibitem[Sch80]{Schmidt80}
W.M. Schmidt.
\newblock {\em Diophantine approximation}, volume 785 of {\em Lect. Notes in
  Math.}
\newblock Springer-Verlag, Berlin--Heidelberg, 1980.

\bibitem[Sch89]{Scheurle89}
J.~Scheurle.
\newblock Chaos in a rapidly forced pendulum equation.
\newblock In {\em Dynamics and Control of Multibody Systems (Brunswick, ME,
  1988)}, volume~97 of {\em Contemp. Math.}, pages 411--419. Amer. Math. Soc.,
  Providence, RI, 1989.

\bibitem[Sim94]{Simo94}
C.~Sim\'o.
\newblock Averaging under fast quasiperiodic forcing.
\newblock In J.~Seimenis, editor, {\em Hamiltonian Mechanics: Integrability and
  Chaotic Behavior (Toru\'n, 1993)}, volume 331 of {\em NATO ASI Ser.~B:
  Phys.}, pages 13--34. Plenum, New York, 1994.

\bibitem[Spi99]{Spinadel99}
V.W.{ }de Spinadel.
\newblock The family of metallic means.
\newblock {\em Vis. Math.}, 1(3):\ approx. 16~pp., 1999.
\newblock \url{http://members.tripod.com/vismath/}.

\bibitem[SV01]{SimoV01}
C.~Sim\'o and C.~Valls.
\newblock A formal approximation of the splitting of separatrices in the
  classical {A}rnold's example of diffusion with two equal parameters.
\newblock {\em Nonlinearity}, 14(6):1707--1760, 2001.

\bibitem[Tre94]{Treschev94}
D.V. Treschev.
\newblock Hyperbolic tori and asymptotic surfaces in {H}amiltonian systems.
\newblock {\em Russian J.~Math. Phys.}, 2(1):93--110, 1994.

\end{thebibliography}

\def\noopsort#1{}

\end{document}